\documentclass[journal]{IEEEtran}
\usepackage{amsmath,amssymb,mathrsfs}
\usepackage{algorithm}
\usepackage{algorithmic}
\usepackage{epsfig}
\usepackage{epstopdf}
\usepackage{graphicx}
\usepackage{color}
\newcommand{\EE}{\mathbb{E}}
\newcommand{\RR}{\mathbb{R}}

\newcommand{\Prob}{\mathbb{P}}
\newcommand{\Min}{\textrm{minimize}}
\newtheorem{lemma}{Lemma}

\newtheorem{theorem}{Theorem}

\newtheorem{definition}{Definition}
\newtheorem{assumption}{Assumption}
\newtheorem{proposition}{Proposition}

\ifCLASSOPTIONcompsoc
  \usepackage[nocompress]{cite}
\else
  \usepackage{cite}
\fi
\ifCLASSINFOpdf
\else
\fi
\begin{document}
%
\title{Decentralized Markov Chain Gradient Descent}

\author{Tao Sun,  and  Dongsheng Li
\IEEEcompsocitemizethanks{
\IEEEcompsocthanksitem T. Sun and D. Li are with the  College of Computer, National University of Defense Technology (nudtsuntao@163.com;dsli@nudt.edu.cn).
\IEEEcompsocthanksitem D. Li is the corresponding author.
}

\thanks{This work  is sponsored in part by National Key R\&D Program of China (2018YFB0204300),  and
the National Natural Science Foundation of China under Grants (61932001 and 61906200).}}

\markboth{IEEE  (submitted \today)}%
{Shell \MakeLowercase{\textit{et al.}}: Bare Demo of IEEEtran.cls for Computer Society Journals}

\IEEEtitleabstractindextext{
\begin{abstract}
Decentralized stochastic gradient method emerges as a promising solution for solving large-scale machine learning problems.
This paper studies the decentralized Markov chain gradient descent (DMGD) algorithm --- a variant of the decentralized stochastic gradient methods where the random samples are taken along the trajectory of a Markov chain. This setting is well-motivated when obtaining independent samples is costly or impossible, which excludes the use of the traditional stochastic gradient algorithms. Specifically, we consider the first- and zeroth-order versions of decentralized Markov chain gradient descent over a connected network, where each node only communicates with its neighbors about intermediate results. The nonergodic convergence and the ergodic convergence rate of the proposed algorithms have been rigorously established, and their critical dependences on the network topology and the mixing time of Markov chain have been highlighted.  The numerical tests further validate the sample efficiency of our algorithm.
\end{abstract}

}

\maketitle

\IEEEdisplaynontitleabstractindextext

%
\IEEEpeerreviewmaketitle

\section{Introduction}
Distributed machine learning is an attractive solution to tackle large-scale learning tasks \cite{dean2012large,li2014scaling}.
In this paper, we consider that $m$ agents represent the nodes of a connected network and collaboratively solve the following stochastic optimization problem
\begin{equation}\label{minex}
    \min_{x\in\mathbb{R}^n}~~~\sum_{i=1}^m ~\EE_{\xi(i)}\big(F(x;\xi(i))\big),
\end{equation}
where $\EE_{\xi(i)}\big(F(x;\xi(i))\big):=\int_{\Pi_i} F(x,\xi(i)) d\Pi_i(\xi(i))$ and $\Xi_i$ is a statistical sample space with probability distribution $\Pi_i$ at node $i$ (we omit the underlying $\sigma$-algebra), and  $F(\cdot;\xi(i)):\mathbb{R}^n\rightarrow \mathbb{R}$ is a closed (possibly nonconvex) function associated with $\xi(i)\in\Xi_i$.
This formulation contains various multi-agent machine learning, reinforcement learning and statistical problems.
We are particularly interested in cases where obtaining an independent and identically distributed (i.i.d.) sample $\xi(i)$ from $\Xi_i$ is very hard or even impossible at every node $i$; see an example in \cite{sun2018markov} where the cost of i.i.d. sampling can be very expensive. In statistics, a common method to overcome the issue of difficult sampling is employing a Markov chain whose stationary distribution is $\Pi_i$.
Therefore, to solve \eqref{minex}, one can still use the parallel implementation of the widely-used method Stochastic Gradient Descent (SGD) \cite{robbins1951stochastic}:
\begin{align}\label{dsgd}
    x^{k+1}=x^k-\gamma_k\sum_{i=1}^m\nabla F(x^k;\xi^k(i)),\quad \xi^k(i) \sim \Pi_i
\end{align}
where $\gamma_k$ and $x^k$ are the stepsize and parameter at iteration $k$, and $\nabla F(x^k;\xi^k(i))$ is a \emph{nearly unbiased} stochastic gradient of $\EE_{\xi(i)}\big(F(x;\xi(i))\big)$ obtained at node $i$.
By using the Markov chain, in order to perform one iteration of \eqref{dsgd}, each node $i$ has to generate a sequence of samples $\xi^1(i),\xi^2(i),\ldots,\xi^T(i)$  and only uses the last one  $\xi^k(i):=\xi^T(i)$.  According to  \cite[Theorem 4.9]{montenegro2006mathematical}, to get a sample that is nearly i.i.d., one needs to simulate the Markov chain for a sufficiently long time; e.g., a large $T$. For this reason, we call the iteration \eqref{dsgd} as SGD-$T$.
In addition, to update the next parameter via \eqref{dsgd}, all the local gradients need to be collected.
Therefore, applying iteration \eqref{dsgd} for \eqref{minex} over the distributed nodes has following two limitations.

\noindent \textbf{Sample inefficiency.} Different from standard stochastic optimization settings, when it is difficult to obtain i.i.d. samples $\xi(i)$ from $\Xi_i$, implementing SGD-$T$ for \eqref{minex} requires regenerating Markov chains at each node per iteration. Nevertheless, this wastes a sizeable amount of variable samples, especially when the Markov chain has a large mixing time.

\noindent \textbf{Communication inefficiency.} The presumption of implementing \eqref{dsgd} is that there is a fusion center (which can be a designated agent) aggregating the local gradients and carrying out the parameter update. However, this incurs a significant amount of synchronization and communication overhead, especially when the network is large and sparse.


 In this context, our goal is to find the near-optimal solution of \eqref{minex} in a sample- and communication-efficient manner.

\subsection{Prior Art}
In this part, we briefly review three lines of related works: decentralized optimization, decentralized stochastic optimization, and Markov chain gradient descent.

\noindent\textbf{Decentralized optimization.} Decentralized algorithms have been originally studied in control and signal processing communities, e.g., calculating the mean of data distributed over multiple sensors \cite{boyd2005gossip,olfati2007consensus,schenato2007distributed,aysal2009broadcast}. The decentralized (sub)gradient descent
 (DGD) algorithms for the finite sum optimization have been studied in \cite{nedic2009distributed,chen2012fast,jakovetic2014fast,matei2011performance,yuan2016convergence}. However, DGD converges to an
inexact solution due to that it actually minimizes an unconstrained penalty function rather than the original one. To fix this, the dual information is leveraged in recent works such as decentralized ADMMs and primal-dual algorithms \cite{chang2015multi,schizas2008consensus,shi2014linear,shi2015extra}.
Although DGD is slower than decentralized ADMMs and primal-dual algorithms in the convex settings,
it is much simpler and therefore easier to extend to the nonconvex, online and delay-tolerant settings \cite{zeng2018nonconvex,hosseini2016online,mcmahan2014delay}.

\noindent\textbf{Decentralized stochastic optimization.}
Generalizing methods for the decentralized deterministic optimization, decentralized SGD (DSGD) has been studied recently. By assuming  a
local Poisson clock for each agent, asynchronous gossip algorithms is proposed by \cite{ram2010asynchronous2}, in which each worker randomly selects part of its neighbors to communicate with. In fact, these algorithms used random communication graphs.
The decentralized algorithm with random communication graph for the constrained problem is introduced by \cite{srivastava2011distributed}; the subgradient version is given by\cite{ram2010distributed}.
In recent works \cite{sirb2016consensus,lan2017communication,lian2017can}, the  theoretical convergence complexity analysis  of  convex and nonconvex DSGD   is provided.
In \cite{sirb2016consensus}, the  authors presented the  complexity analysis for  a
stochastic  decentralized   algorithm.  In \cite{lan2017communication}, the authors design a kind of stochastic  decentralized   algorithm  by recruiting the dual information, and provide the related computational complexity. In latter paper \cite{lian2017can}, the authors show the  speedup when the number of nodes is increased. And in paper \cite{lian2018a}, the authors proposed  the asynchronous decentralized stochastic gradient descent and presented the theoretical and numerical results.

\noindent\textbf{Markov chain gradient descent.}
While i.i.d. samples are not always available in stochastic optimization, recent focus has been on the analysis of stochastic algorithms following a single trajectory of the Markov chain or other general ergodic processes.
The key challenge of analyzing MGD is to deal with the biased expectation of gradients. The ergodic convergence results have been reported in  \cite{johansson2009randomized,johansson2007simple}. Specifically, \cite{johansson2009randomized,johansson2007simple} study the conditional expectation with a sufficiently large delay which is sufficiently close to the gradient (but still different). The authors of \cite{ram2009incremental}
proved the almost sure convergence under the diminishing stepsizes  $\gamma_k={1}/{k^q}$, ${2}/{3}< q\leq 1$.  In \cite{duchi2012ergodic}, the authors improved convergence results with larger stepsizes $\gamma_k={1}/{\sqrt{k}}$ in the sense of ergodic convergence.
In all these works, the Markov chain is required to be reversible, and the functions have to be convex. In a very recent paper \cite{sun2018markov}, the non-ergodic convergence of MGD has been shown in the nonconvex case with non-reversible Markov chain, but the algorithm needs to be implemented in a centralized fashion.

\subsection{Our contributions}
The main theme of this paper is the development of the first- and zeroth-order version of decentralized Markov chain gradient descent (DMGD) algorithms, and their performance analysis. In contrast to the well-known DSGD, DMGD leverages Markov chain sampling rather than uniform random sampling, which gains sample efficiency. For first-order DMGD, each node uses a Markov chain trajectory to sample a gradient and then communicates with its neighbors to update the local variables. To further account for the case where stochastic gradient is not readily available, we develop the zeroth-order variant of DMGD where knowing the gradients is not necessary, and only the point-wise function values are needed.

 We establish the non-ergodic convergence of first- and zeroth-order DMGD and their ergodic convergence rates.  The results show that DMGD --- the first-order DMGD converges at the same order as the centralized MGD, and the zeroth-order DMGD converges at the same order as DMGD. Some novel results are are developed based on new techniques and approaches developed in this paper. To get the stronger results in general cases, we used the varying mixing time rather than fixed ones.
It is worth mentioning that our theoretical results can directly derive novel results for DSGD if the Markov chain trajectory sampling reduces to the uniform random sampling case. The numerical results are presented to demonstrate that DMGD performs better than DSGD in terms of sample efficiency.

Although this paper's proofs also need to use the delay expectation techniques used in \cite{sun2018markov}, the significant difference remains.
This is because the DMGD introduced in this paper  has no objective function to minimize.
The proof of DMGD is built on the average of all nodes' parameters, which is also quite different from MCGD. To this end, several techniques are developed to characterize the difference between the average
and nodes' parameters under Markov chain sampling.
Compared with  \cite{sun2018markov}, this paper also investigates the zeroth-order case that has never been studied in MCGD.

\vspace{0.2cm}

\noindent\textbf{Notation}:
Let $\lambda_{i}(\cdot)$ denote the $i$-th eigenvalue of a matrix. Let ${\bf x}(i)\in \mathbb{R}^N$  denote the local copy of $x$ at node $i$.
For a matrix $A=(a_{i,j})_{m\times n}$,   $\|A\|:=\sqrt{\sum_{i=1}^m\sum_{j=1}^n a_{i,j}^2}$ and $\|A\|_{\infty}:=\max_{i,j}|a_{i,j}|$. For a positive semidefinite matrix $B$, $\|A\|_{B}:=\|B^{\frac{1}{2}}A\|$.

\section{Decentralized Markov Chain Gradient Descent}
\subsection{Preliminaries}
We first consider the discrete case of our problem \eqref{minex}, that is,  all the distributions  ($\Pi_1,\Pi_2,\ldots,\Pi_m$) are supported on a set of $M$ points\footnote{For notational brevity, we assume the same cardinal number for different distribution support sets.}, $y^{1,j},\dots,y^{M,j}$ (for $\Pi_j$). We define the functions as
$f_{j}^i(x) := M\cdot\mathrm{Prob}(\xi=y^{i,j})\cdot F(x,y^{i,j})$, and thus problem \eqref{minex} becomes the following finite-sum formulation
\begin{equation}\label{summodel}
    \Min_{x\in\mathbb{R}^{n}}~~f(x):=\frac{1}{m}\sum_{j=1}^m f_j(x),
\end{equation}
where $f_j(x):=\frac{1}{M} \sum_{i=1}^{M} f_{j}^i(x)$  is the loss function in the $j$th node. The consensus optimization is the main workhorse to distribute the training algorithm in machine learning tasks. There are just a few non-consensus methods; as far as we know, it seems only block coordinate gradient descent\cite{wright2015coordinate,tseng2001convergence,richtarik2016distributed}. Mathematically, it applies the coordinate descent to minimize function $f$ or its dual problem. Thus, we can get a centralized system whose work computes the gradient for one block for algorithmic distribution. However, block coordinate gradient descent has still not been widely used for distributed training tasks due to the vast data may cause random access memory overflow even for one block; on the other hand, the large batch size implies indecent generalization property \cite{keskar2016large} that further limits the use of block coordinate gradient descent. Hence, this paper still follows the consensus optimization routine.


Denote $(j_{i,k})_{k\geq 0}\subseteq\{1,2,\ldots,M\}$ as the trajectory of the Markov chain in the $i$th node  and $k$th iteration. We use a connected graph
$\mathcal{G} = (\mathcal{V}, \mathcal{E})$ with vertex set $\mathcal{V}=\{1,...,M\}$ and edge set $\mathcal{E}\subseteq \mathcal{V}\times \mathcal{V}$.
Any edge $(i, l)\in\mathcal{E}$ represents a communication link between nodes  $i$
and $l$.
 And more, let
 $$
\mathbf{j}_k:=\begin{bmatrix}
    j_{1,k} \\
    j_{2,k} \\
    \vdots \\
    j_{m,k} \\
\end{bmatrix},~
{\bf x}:=\begin{bmatrix}
    {\bf x}(1)^{\top} \\
    {\bf x}(2)^{\top} \\
    \vdots \\
    {\bf x}(m)^{\top} \\
\end{bmatrix},~
{\bf u}^k:=\begin{bmatrix}
    \nabla f_1^{j_{1,k}}({\bf x}^k(1)) \\
    \nabla f_2^{j_{2,k}}({\bf x}^k(2)) \\
    \vdots \\
   \nabla f_m^{j_{m,k}}({\bf x}^k(m))\\
\end{bmatrix}.$$


\textbf{Mixing matrix}:
The mixing matrix is frequently used in decentralized optimization. In many cases, it can be designed by the users according to the given graph.
Formally, it is defined as follows.
\begin{definition}
\label{def:MixMat}
The mixing matrix $W = [w_{ij}] \in \mathbb{R}^{m\times m}$ is assumed to have the following properties:\\
(1) {\em (Graph)} If $i\neq j$ and $(i,j) \notin {\cal E}$, then $w_{ij} =0$, otherwise, $w_{ij} >0$;\\
(2) ({\em Symmetry}) $W = W^{\top}$;\\
(3) ({\em Null space property}) $\mathrm{null} \{I-W\} = \mathrm{span}\{\bf 1\}$;\\
(4) ({\em Spectral property}) $I \succeq W \succ -I.$
\end{definition}
With the symmetricity of $W$, its eigenvalues are real and can be sorted in the nonincreasing order.
Thus, let  $\lambda_i(W)$ denote the $i$th largest eigenvalue of $W$; then, it holds that
$\lambda_1(W)=1>\lambda_2(W) \geq \cdots \geq \lambda_m(W)>-1.$
%

\textbf{Markov chain}:
We recall several definitions, properties, and existing results of the
finite-state time-homogeneous Markov chain, which will be used in the proposed algorithms.
\begin{definition}
Let $H$ be an $n\times n$-matrix with real-valued elements.
A stochastic process $X_1,X_2,...$ in a finite state space $\{1,2,\ldots,n\}$ is called a time-homogeneous Markov chain with transition matrix $H$ if, for $k\in \mathbb{N}$, $i,j\in \{1,2,\ldots,n\}$, and  $i_0,i_1,\ldots,i_{k-1}\in \{1,2,\ldots,n\}$, we have
\begin{align*}
    &\mathbb{P}(X_{k+1}=j\mid X_0=i_0,X_1=i_1,\ldots,X_k=i)\\
    &\qquad=\mathbb{P}(X_{k+1}=j\mid X_k=i)=H_{i,j}.
\end{align*}
\end{definition}
Denote  the probability distribution of $X_k$  as the non-negative row vector $\pi^k=(\pi^k_1,\pi^k_2,\ldots,\pi^k_n)$, i.e., $\Prob(X_k=j)=\pi_j^k$ and $\pi$ satisfies  $\sum_{i=1}^n \pi^k_i=1.$
For the time-homogeneous Markov chain, it holds $\pi^k=\pi^{k-1} H$ and
$
    \pi^k=\pi^{k-1} H=\cdots=\pi^0 H^{k},
$
for $k\in \mathbb{N}$, where $H^{k}$ denotes the $k$th power of $H$.

A Markov chain is irreducible if, for any $i,j\in \{1,2,\ldots,n\}$, there exists $k$ such that $(H^k)_{i,j}>0$.
State $i\in\{1,2,\ldots,n\}$ is said to have a period $d$ if $H^k_{i,i} = 0$ whenever $k$ is \emph{not} a
multiple of $d$ and $d$ is the greatest integer with
this property. If $d=1$, then we say state $i$ is aperiodic. If  every state is aperiodic, the Markov chain is said to be aperiodic.
Any time-homogeneous, irreducible, and aperiodic
Markov chain has a stationary distribution
$\pi^*=\lim_k \pi^k
=[\pi^*_1,\pi^*_2,\ldots,\pi^*_n]$
 with $\sum_{i=1}^n \pi^*_i=1$ and $\min_i\{\pi^*_i\}>0$, and $\pi^*= \pi^* H$. It also holds that
\begin{equation}\label{convermatrix}
    \lim_{k}H^{k}=\begin{bmatrix}(\pi^*)^{\top},(\pi^*)^{\top},\ldots,(\pi^*)^{\top}\end{bmatrix}^{\top}
    =:\Pi^*\in \mathbb{R}^{n\times n}.
\end{equation}
The largest eigenvalue of $H$ is 1, and the corresponding  eigenvector is $\pi^*$.

Mixing time is an important concept which describes how long a Markov chain evolves until its current state has a distribution very close to its stationary distribution. The literature    studies about various kinds of mixing times, whose majority, however, is about
 reversible Markov chains (i.e.,  $\pi_i H_{i,j} = \pi_j H_{j,i}$).  With basic matrix analysis, the  mixing time introduced in \cite{sun2018markov} provides a direct relationship between $k$ and the deviation of the distribution of the current state from the stationary distribution (Lemma \ref{lemc} in the Appendix).

 \subsection{Algorithmic Development of DMGD}
The local scheme of DMGD at the $i$th node is
\begin{align}\label{localscheme}
{\bf x}^{k+1}(i)=\sum_{l\in \mathcal{N}(i)} w_{i,l} {\bf x}^{k}(l)-\gamma_k\nabla f^i_{j_{i,k}}({\bf x}^k(i)).
\end{align}
In each iteration, each node  calculates the local gradient on the Markov chain trajectory $(j_{i,k})_{k\geq 0}$, and then communicates with its neighbors $\mathcal{N}(i)$ with  a weighted average $\sum_{l\in \mathcal{N}(i)} w_{i,l} {\bf x}^{k}(l)$ to update the iteration. Here, $w_{i,l}$ is the $(i,l)$-element of the mixing matrix. It is easy to see that if the  Markov chain trajectory is the uniform sampling, \eqref{localscheme} then reduces to the DSGD.
The parameter $\gamma_k$ will go to zero. The settings are to guarantee  the convergence of the algorithm; the diminishing stepsizes are used to reduce the variances cost by the gradients samplings.
We present the comparison of DSGD-$T$ (Algorithm \ref{alg1}) and DMGD (Algorithm \ref{alg2}).
Our work is the first to distribute the consensus optimization under the Markov sampling assumption in the decentralized case. In our algorithm (8), if $\mathcal{N}(i)=\{i\}$, it is then the centralized Markov chain gradient descent. Thus, our method contains the previous distributed method for Markov chain samplings.

\begin{algorithm}
\caption{DSGD-$T$}
\begin{algorithmic}\label{alg1}
\REQUIRE   parameters  $(\gamma_k)_{k\geq 0}$\\
\textbf{Initialization}: ${\bf x}(1)={\bf x}(1)=\ldots={\bf x}(m)=x^0$\\
\textbf{for}~$k=1,2,\ldots$ \\
~~~\textbf{for}~$i=1,2,\ldots,m$ \\
 \STATE  1. \textbf{Resample} a Markov chain $j_0(i),\ldots,j_{T}(i)$\\
 \STATE 2. Collect ${\bf x}^k(i)$ from its neighbors \\
 \STATE 3. Update ${\bf x}^{k+1}(i)$ via \eqref{localscheme} with $j_{i,k}\leftarrow j_{T}(i)$\\
~~~\textbf{end for}\\
\textbf{end for}\\
\end{algorithmic}
\end{algorithm}

\begin{algorithm}
\caption{DMGD}
\begin{algorithmic}\label{alg2}
\REQUIRE   parameters  $(\gamma_k)_{k\geq 0}$\\
\textbf{Initialization}: ${\bf x}(1)={\bf x}(1)=\ldots={\bf x}(m)=x^0$\\
\textbf{for}~$k=1,2,\ldots$ \\
~~~\textbf{for}~$i=1,2,\ldots,m$ \\
 \STATE  1. \textbf{Sample} $j_{i,k}$ via a   Markov chain\\
 \STATE 2. Collect ${\bf x}^k(i)$ from its neighbors \\
 \STATE 3. Update ${\bf x}^{k+1}(i)$ via \eqref{localscheme} with $j_{i,k}$\\
~~~\textbf{end for}\\
\textbf{end for}\\
\end{algorithmic}
\end{algorithm}

The global scheme can be described as
\begin{equation}\label{dmgd2}
	    {\bf x}^{k+1}= W{\bf x}^k-\gamma_k {\bf u}^k
\end{equation}
where  ${\bf u}^k$ has been given before.
This iterative formulation   can help us to understand the convergence of the algorithm.
Suppose that the Markov chains all reduce to the uniform sampling.
By defining the $\sigma$-algebra as
$\chi^k:=\sigma({\bf x}^0,{\bf x}^1,\ldots,{\bf x}^k),$
we can see  that $\EE({\bf u}^k\mid\chi^k)= \sum_{i=1}^{m} \nabla f_i[{\bf x}^k(i)]$ in this condition;  that means DMGD
actually admits the DSGD  with converge guarantees   in the perspective of expectation.  For general Markov chains, the analysis is much more complicated for  the biased conditional expectation.

 \subsection{Key Challenge of Analyzing DMGD}
Markov chain sampling   is neither cyclic nor i.i.d. stochastic. For any large $K$, it is still possible that a sample is never visited during some $k+1,...,k+K$ iterations. For a fixed node $i$, unless the graph $\mathcal{G}_i$ is a  complete graph (i.e., all elements are directly connected), there are elements $l,h$ \emph{without} an edge connecting them, i.e., $(l_i,h_i)\not\in \mathcal{E}_i$. Hence, given $j_{i,k-1} = l_i$, it is \emph{impossible} to have $j_{i,k}=h_i$. So, no matter how one selects the sampling probability   and stepsize $\gamma_k$, we generally  \emph{cannot} have
$     \EE(\gamma_k{\bf u}^k\mid {\bf j_{k-1}}=\mathbf{l})= C (\sum_{i=1}^{m} \nabla f_i[{\bf x}^k(i)])$
for any constant $C$, where $\mathbf{l}=(l_1,\ldots,l_m)$. This fact, unfortunately, breaks down all the existing analyses of stochastic decentralized optimization since they all need a non-vanishing probability for every sample in each node can be selected.

 \section{Convergence analysis of DMGD}
In this part, we present the theoretical results of DMGD with finite-state Markov chains.
Our analysis builds on the following assumptions.
 \begin{assumption}\label{ass0}
Function $f_j$ is lower bounded, that is
$f_j(x)> -\infty,\,\forall j, \forall x$.
\end{assumption}
\begin{assumption}\label{ass1}
The gradient of $f^i_j$ is uniformly bounded, that is, there exists a constant $B>0$ such that
$
    \|\nabla f^i_j(x)\|\leq B, \forall i\in\{1,2,\ldots,M\},~\forall j\in\{1,2,\ldots,m\}
$.
\end{assumption}
\begin{assumption}\label{ass2}
The gradient of $f^i_j$ is Lipschitz continuous with $L_j^i$, i.e.,
$
    \|\nabla f^i_j(x)-\nabla f^i_j(x)\|\leq L^i_j\|x-y\|, \forall i\in\{1,2,\ldots,M\},~\forall j\in\{1,2,\ldots,m\}
$. We denote that $L:=\max_{1\leq i\leq m, 1\leq j\leq M}\{L^i_j\}$.
\end{assumption}

\begin{assumption}\label{ass:mc}{The Markov chains in all nodes are time-homogeneous, irreducible, and aperiodic, which  have the same transition matrix $H$ and the same stationary distribution.}\footnote{We require all nodes to employ the Markov chain with same transition matrix $H$. This setting is for the convenience of presentations in the proofs and can be modified as different Markov chains.}
\end{assumption}

Following the routines in stochastic decentralized community, the convergence of the algorithm is described by
$$\overline{x^k}:=\frac{1}{m}{\sum_{j=1}^m {\bf x}^k(j)}.$$
\begin{theorem}\label{th2}
Suppose Assumptions \ref{ass0}-\ref{ass:mc} hold and the stepsizes are selected as
\begin{equation}\label{stepsizes}
\gamma_k=\frac{1}{(k+1)^{\theta}},~~~\frac{1}{2}<\theta<1.
\end{equation}
For $({\bf x}^k)_{k\geq 0}$ generated by DMGD, we have the following nonergodic convergence result
\begin{equation}\label{th2-r1}
    \lim_{k}\EE\|\nabla f(\overline{x^k})\|=0.
\end{equation}
And the ergodic convergence rate is
\begin{align}\label{th2-r2}
    &\min_{1\leq i\leq k}\EE\|\nabla f(\overline{x^k})\|={\cal O}\left(\frac{\frac{1}{\ln(1/\lambda(H))}\cdot\frac{1}{1-\lambda_2(W)}}{k^{\frac{1-\theta}{2}}}\right),
\end{align}
where $\lambda(H):=\frac{\max\{|\lambda_2(H)|,|\lambda_{\min}(H)|\}+1}{2}\in [0,1)$, and $\lambda_2(H)$ and $\lambda_{\min}(H)$ denote the second and smallest eigenvalue of $H$, respectively.
\end{theorem}
In Theorem \ref{th2}, the functions are not necessary to be convex.
In fact, it  is more difficult to prove \eqref{th2-r1} than to prove \eqref{th2-r2}. The descent on a Lyapunov function and the Schwarz inequality can directly derive \eqref{th2-r1}, while \eqref{th2-r2} requires a technical lemma, which first given in \cite{zeng2018nonconvex} and generalized by \cite{sun2018markov}. An extreme case is that $m=1$ and $W=1$; DMGD will reduce to the classical MGD. But Theorem \ref{th2} cannot  cover  the existing  convergence results of MGD. In \cite{sun2018markov}, the authors estimated the convergence  of MGD with the stepsizes requirements
\begin{equation}\label{stepsizes2}
\sum_{k=1}^{+\infty}\gamma_k=+\infty,~~~\sum_{k=1}^{+\infty}\ln^2 k\cdot\gamma_k^2<+\infty.
\end{equation}
The stepsize \eqref{stepsizes} can satisfy \eqref{stepsizes2}  but not vice versa.

The convergence results in  Theorem \ref{th2} can be extended to general Markov chains with extra assumptions given in existing works.
In \cite{ram2009incremental}, the time non-homogeneous Markov chain but with extra assumptions (Assumptions 4 and 5, in Section 4 of \cite{ram2009incremental}) is proposed. These two assumptions involve with many details; several majors are doubly stochastic, uniformly bounded away from zero, diagonals of the transition matrices are  positive, and strong connections of some edges. In paper \cite{duchi2012ergodic}, the authors use more  general Markov chain but also with an assumption (Assumption C, in Section 2 of \cite{duchi2012ergodic}), which can be satisfied by finite-state time-homogeneous Markov chain.

 To reach the error as $\EE\|\nabla f(\overline{{\bf x}^k})\|^2\leq \epsilon$, we need the learning rate as
$\gamma=\Theta(\epsilon)$ and $K=\widetilde{O}(\frac{1}{\epsilon^2})$. This result is almost the same as the diminishing learning rate case. However, the constant learning rate case cannot prove the non-ergodic convergence result (11). Thus, this paper uses diminishing ones.




\section{Zeroth-order DMGD}
This section presents the zero-order version  of DMGD with two-points feedback strategy \cite{duchi2015optimal,ghadimi2013stochastic,agarwal2010optimal,shamir2017optimal}. This paper employs the method given in \cite{duchi2015optimal,ghadimi2013stochastic}. Specifically, it uses the estimator of the gradient of $f$ by  querying at ${\bf x}+\delta {\bf h}$ and ${\bf x}$ with returning
$
	\frac{n(f({\bf x}+\delta {\bf h})-f({\bf x}))}{\delta}{\bf h}
$
where ${\bf h}$ is a random unit vector and $\delta>0$ is a small parameter. In the zeroth-DMGD, we use the two-points feedback to replace the local gradients and obtain following iteration at the $i$th node
\begin{align}\label{localschemezero}
&{\bf x}^{k+1}(i)=\sum_{l\in \mathcal{N}(i)}\! w_{i,l} {\bf x}^{k}(l)\nonumber\\
&-\frac{\gamma_k n(f^i_{j_{i,k}}({\bf x}^k(i)+\delta_k {\bf h}^{i,k})-f^i_{j_{i,k}}({\bf x}^k(i)))}{\delta_k}{\bf h}^{i,k},
\end{align}
where $(j_{i,k})_{k\geq 0}$ still denotes the    Markov chain trajectory in the $i$th node, and ${\bf h}^{i,k}$ is uniformly sampled for the unite sphere in $\mathbb{R}^n$, and $\delta_k$ is the parameter in the $k$th iteration.  In this algorithm, we just use the function values information rather than the gradients. Thus it is called zeroth-order scheme. We can present the following convergence result of the zeroth-order DMGD.

\begin{algorithm}
\caption{Zeroth-DMGD}
\begin{algorithmic}\label{algo:DMGD2}
\REQUIRE   parameters  $(\gamma_k)_{k\geq 0}$\\
\textbf{Initialization}: ${\bf x}(1)={\bf x}(1)=\ldots={\bf x}(m)=x^0$\\
\textbf{for}~$k=1,2,\ldots$ \\
~~~\textbf{for}~$i=1,2,\ldots,m$ \\
 \STATE  1. \textbf{Sample} $j_{i,k}$ via a   Markov chain\\
 \STATE 2. Collect ${\bf x}^k(i)$ from its neighbors \\
 \STATE 3. Node $i$ samples ${\bf h}^{i,k}$ from the surface of a unit sphere and updates ${\bf x}^{k+1}(i)$ via \eqref{localschemezero}\\
~~~\textbf{end for}\\
\textbf{end for}\\
\end{algorithmic}
\end{algorithm}


\begin{theorem}\label{thzero}
Suppose Assumptions \ref{ass0}, \ref{ass1}, \ref{ass2} and \ref{ass:mc} hold, and ${\bf h}^{i,k}$ is a random unit vector in the $k$th iteration in $i$th node. If the stepsizes are selected as
\begin{equation}
\gamma_k=\frac{1}{(k+1)^{\theta}}\,,\,~~~\frac{1}{2}<\theta<1, ~~~\sum_{k}\gamma_k\delta_k<+\infty
\end{equation}
then for $({\bf x}^k)_{k\geq 0}$ generated by zeroth-order DMGD, we have the nonergodic convergence as
\begin{equation}\label{thzero-r1}
    \lim_{k}\EE\|\nabla f(\overline{x^k})\|=0
\end{equation}
and the ergodic convergence rate is
\begin{align}\label{thzero-r2}
    &\min_{1\leq i\leq k}\EE\|\nabla f( \overline{x^k})\|={\cal O}\left(\frac{\frac{1}{\ln(1/\lambda(H))}\cdot\frac{1}{1-\lambda_2(W)}+n^{\frac{3}{2}}}{k^{\frac{1-\theta}{2}}}\right).
\end{align}
\end{theorem}
Compared with the first-order DMGD in Theorem \ref{th2}, a constant factor $n^{\frac{3}{2}}$ degrades the convergence rates. Such difference comes from the two-points estimation errors of gradient which is dimension-dependent. This result indicates that in low-dimension case, the zeroth-order version can work well as DMGD,  as could be expected; but for high-dimension case, the  speed might be slowed down.

\section{Analysis on continuous  state space}
In this part, we consider the case that $\Pi_1,\Pi_2,\ldots,\Pi_m$  are  continuums, and turn back to the problem \eqref{minex}.
Time-homogeneous and reversible infinite-state Markov chains are considered in this case. With the results in [Theorem 4.9, \cite{montenegro2006mathematical}], the mixing time in this case still enjoys geometric decrease like \eqref{core1}. Mathematically, this fact can be presented as
\begin{equation}\label{core2}
    \|\delta^k\|_{\infty}\leq C\cdot\lambda^k,~~\textrm{as}~~k\geq 0,
\end{equation}
where $\delta^k$ still denotes the deviation matrix $\Pi^*-H^{k}$. Here $C$ and $\lambda$ are constants determined by the Markov chain.  Here, we use notation $C$ and $\lambda$   to give the difference to $C_{P}$ and $\lambda(P)$ in Lemma \ref{lemc}.

Let $\xi^0(i),\xi^1(i),\ldots$ be the Markov chain trajectory  in $i$th node. We first present the   local scheme:
\begin{align}\label{localscheme-online}
&{\bf x}^{k+1}(i)=\sum_{l\in \mathcal{N}(i)} w_{i,l} {\bf x}^{k}(l)-\gamma_k\nabla F({\bf x}^k(i);\xi^k(i)).
\end{align}
By defining $${\bf d}^k:=\begin{bmatrix}\nabla F({\bf x}^k(1);\xi^k(1));\ldots;\nabla F({\bf x}^k(m);\xi^k(m))\end{bmatrix}^{\top},$$ the global scheme is then of the following form
\begin{align}\label{glo-online}
    {\bf x}^{k+1}=W{\bf x}^k- \gamma_k {\bf d}^k.
\end{align}
The convergence is proved for a possibly nonconvex objective function $F(\cdot;\xi)$ and time-homogeneous and reversible chains, which obey the following assumption.
\begin{assumption}\label{assf}
For any $\xi\in\Xi$, it holds that
(1) $\| \nabla F(x;\xi)-\nabla F(y;\xi)\|\leq L\|x-y\|$,  $\forall x,y\in \mathbb{R}^n$;
(2) $\sup_{x\in \mathbb{R}^n,\xi\in\Xi}\{\|\nabla F(x;\xi)\|\}<+\infty$;
(3) $\EE_{\xi}\left(\nabla F(x;\xi)\right)=\nabla \left(\EE_{\xi} F(x;\xi)\right)$, $\forall x\in \mathbb{R}^n$;
(4) $\inf_{x\in \mathbb{R}^n}\left(\EE_{\xi_i} F(x;\xi_i)\right)>-\infty$, $i=1,2,\ldots,m$;
(5) The stationary distribution of the Markov chain in the $i$th node is right $\Pi_i$.
\end{assumption}
Denote a function as
\begin{equation}\label{online-lya}
    \mathcal{F}({\bf x}):=\frac{1}{m}\sum_{i=1}^m \EE_{\xi_i}\big(F({\bf x}(i);\xi_i)\big).
\end{equation}
The convergence results are described by this function.
\begin{proposition}\label{pro-online}
Let Assumption \ref{assf} hold and $({\bf x}^k)_{k\geq 0}$ denote the iterates generated by \eqref{localscheme-online}.
If the stepsizes are selected as \eqref{stepsizes}, we have the following convergence result
\begin{equation}
    \lim_{k}\EE\|\nabla \mathcal{F}( \overline{x^k})\|=0.
\end{equation}
And the ergodic convergence rate is
\begin{equation}
    \min_{1\leq i\leq k}\EE\|\nabla \mathcal{F}( \overline{x^k})\|=O\left(\frac{ \frac{1}{1-\lambda_2(W)}\cdot\frac{1}{\ln(1/\lambda)}}{k^{\frac{1-\theta}{2}}}\right).
\end{equation}
\end{proposition}
Unlike  Theorem \ref{th2}, the Markov chain assumptions cannot be weakened, i.e., the Markov chains must be time-homogeneous and reversible in Proposition \ref{pro-online}. Another difference is that the stationary distributions  $\Pi_1,\Pi_2,\ldots,\Pi_m$ are not necessary to be uniform.
In fact, \eqref{localscheme-online} can be extended to the case where node $i$ stores different functions $F_i(\cdot;\cdot)$; and $F_i(\cdot;\cdot), i=1,2,\ldots,m$ all satisfy Assumption \ref{assf}.

\section{Numerical tests}
\begin{figure*}[!htb]
 \vspace{-0.2cm}
    \centering
    \includegraphics[width=0.3\textwidth]{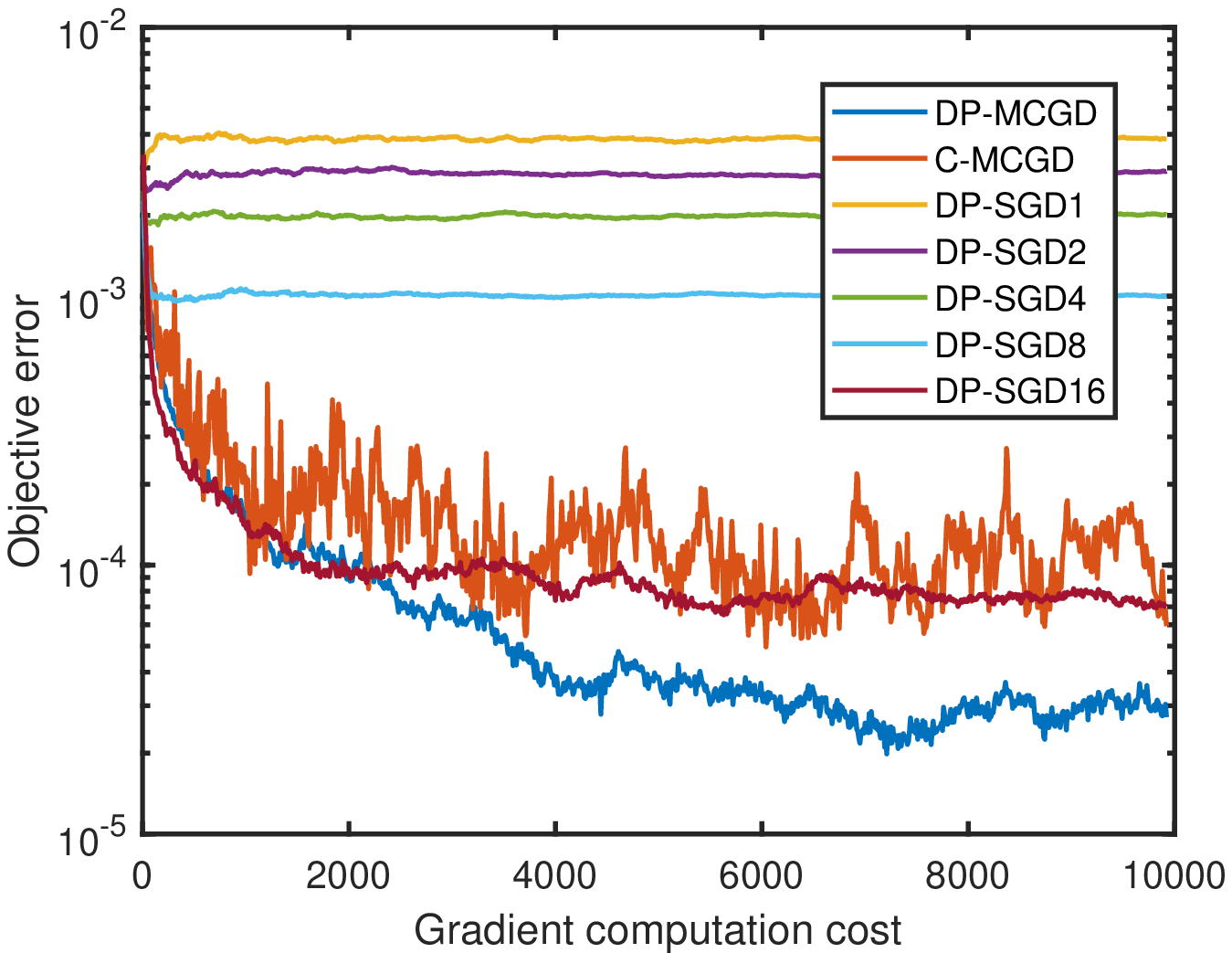}\includegraphics[width=0.3\textwidth]{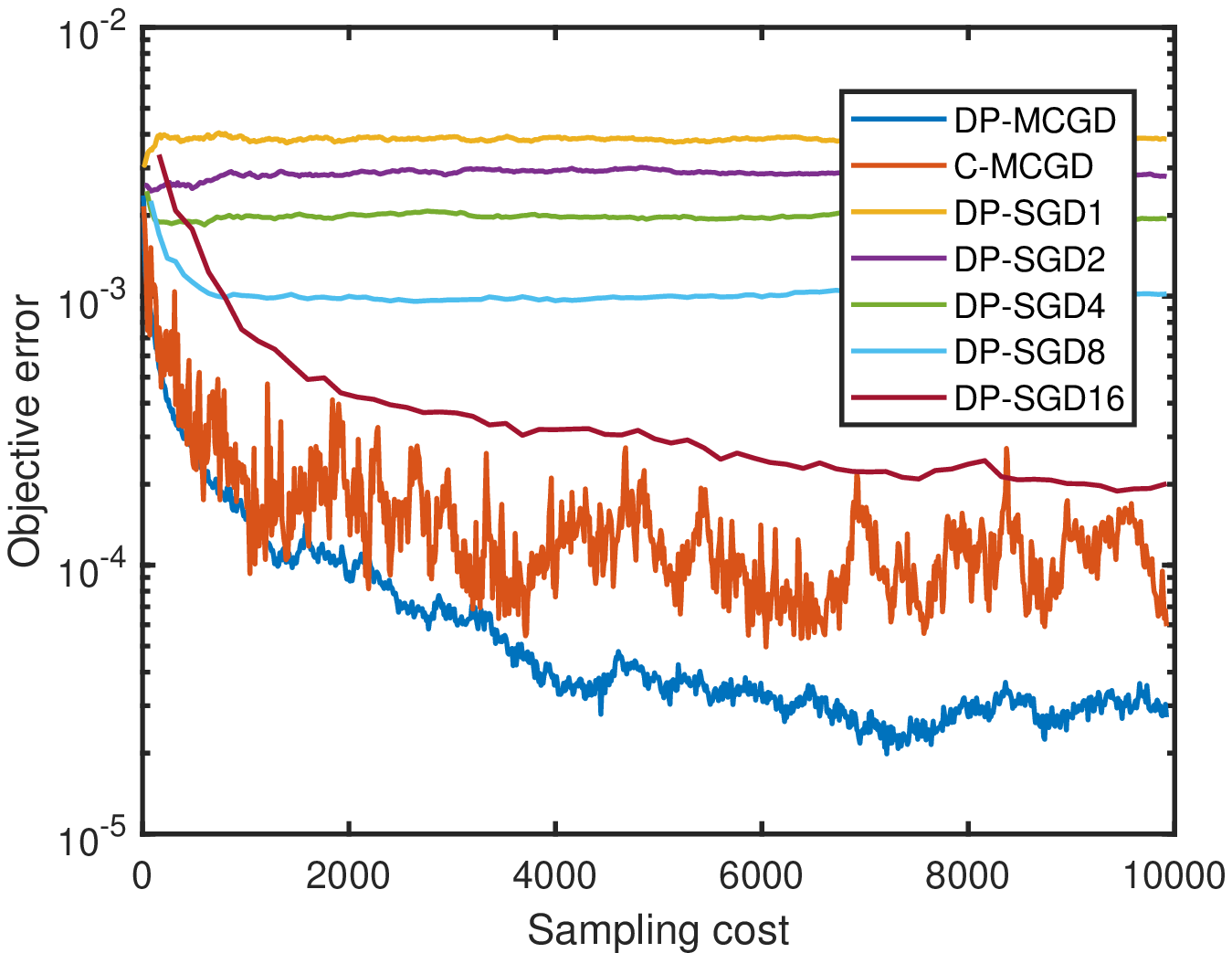}\includegraphics[width=0.3\textwidth]{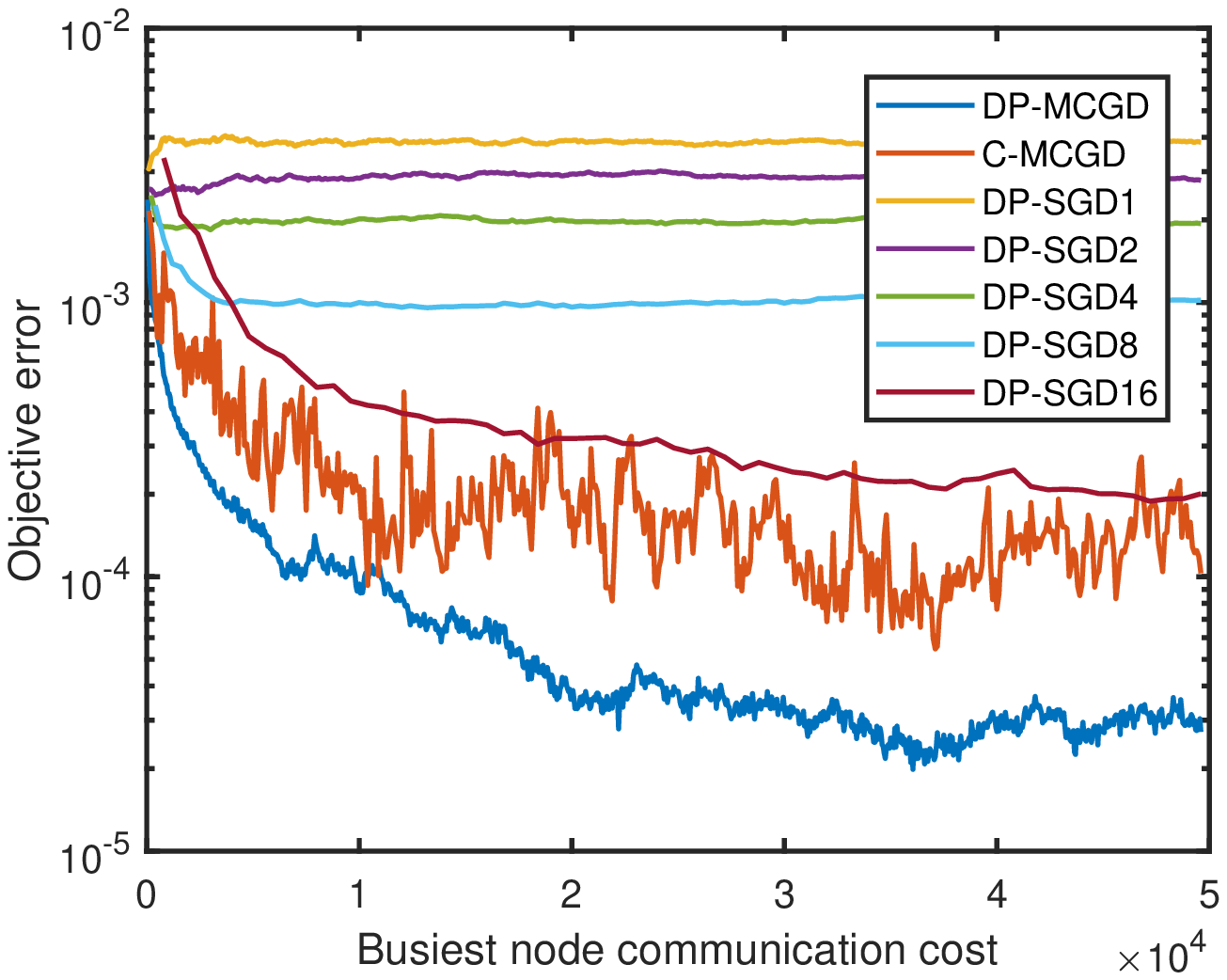}\\
    \vspace{0.6cm}
     \includegraphics[width=0.3\textwidth]{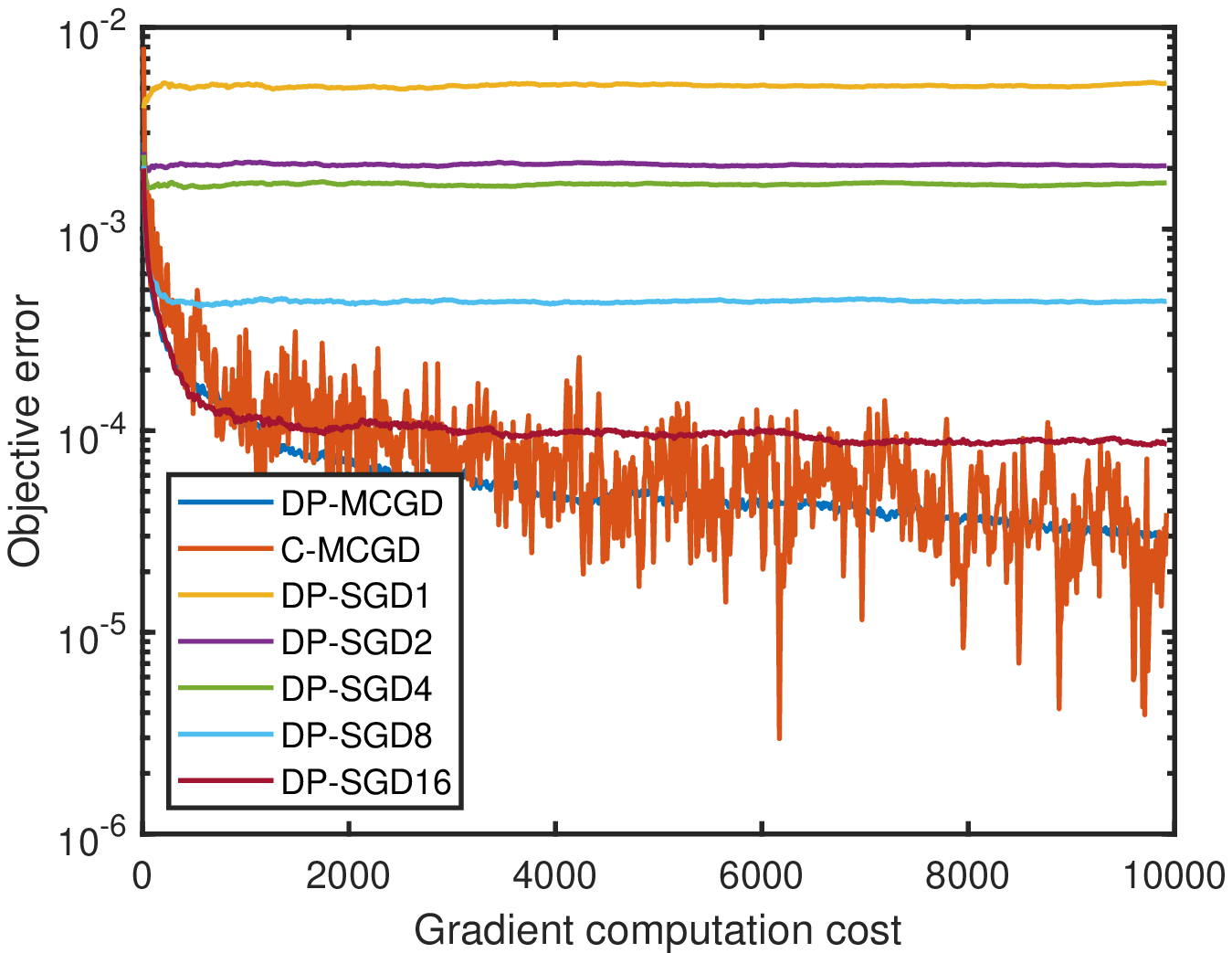}\includegraphics[width=0.3\textwidth]{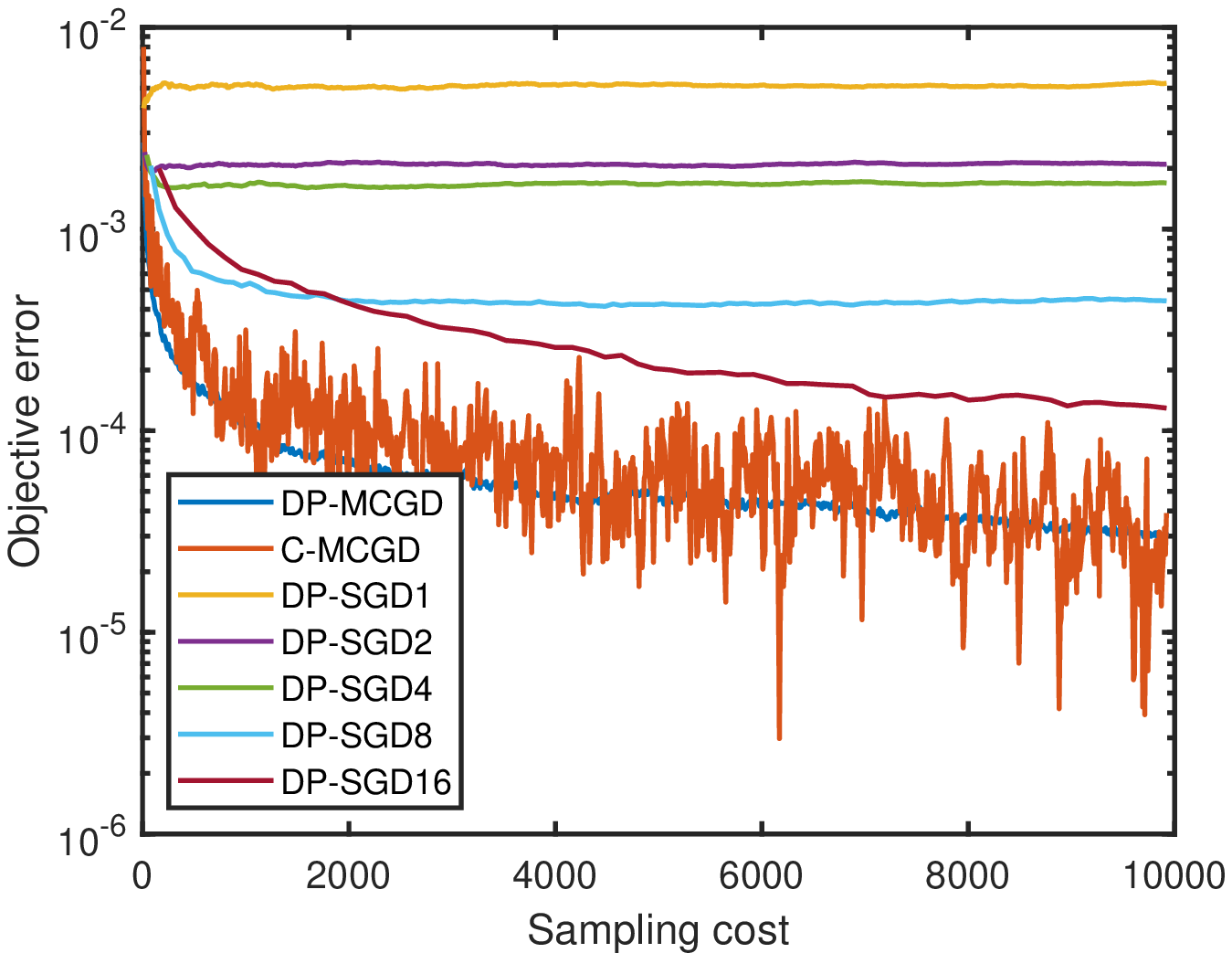}\includegraphics[width=0.3\textwidth]{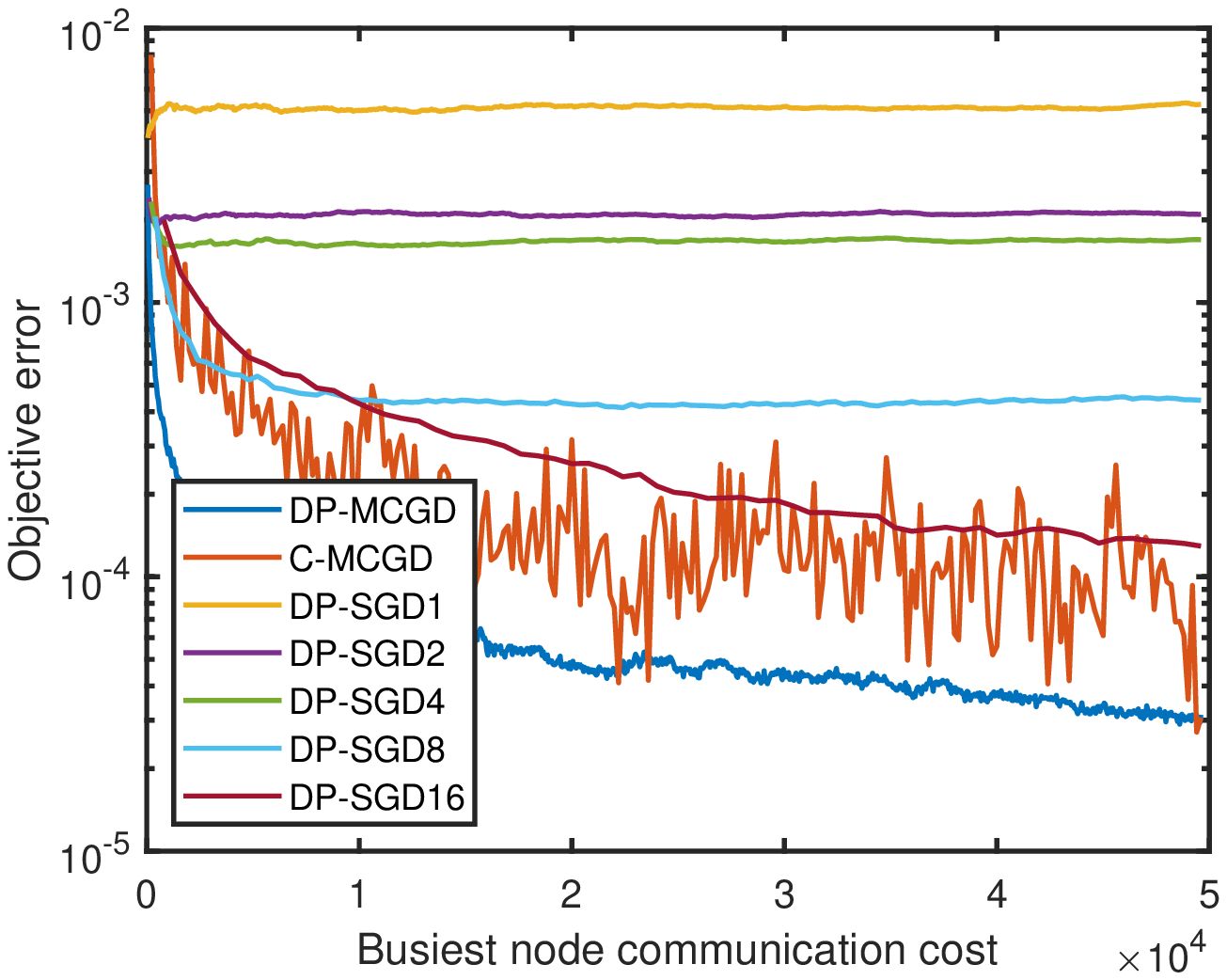}\\
    \caption{Comparisons of DMGD, and centralized MCGD, and DSGD-$T$ for $T=1,2,4,8,16$. The objective function error for (upper row) $m=10,n=50$, and (lower row) $m=20,n=100$.}
    \vspace{-0.2cm}
\end{figure*}
In this section, we compare the performance of our algorithm with the decentralized parallel SGD (DSGD) on an autoregressive model, which closely resembles the first experiment in \cite{duchi2012ergodic}. Assume that there are $m$ autoregressive processes distributed on a graph of $m$ nodes. We attempt to recover a consensus vector $u$ from the multiple processes. On each node $j$, set matrix $A^j$ as a subdiagonal matrix with random entries $A_{l,l-1}(i)\overset{\text{i.i.d}}{\sim} \mathcal{U}[0.8,0.99]$. Randomly sample a vector $u\in\mathbb{R}^n$, with the unit 2-norm. In our experiments, we tested with $m=10, n=50$ and $m=20, n=100$. The data $({\bf\xi}_t^1(i),{\bf\xi}_t^2(i))_{t=1}^{\infty}$ are generated by the following auto regressive process:
\begin{subequations}
	\begin{align*}
    &\xi_t^1(i) = A(i)\xi_{t-1}^1(i)+e_1W^t,~ W^t\stackrel{\text{i.i.d}}{\sim} N(0,1)~~~\,~~\\
    &\bar{\xi}^2_t(i) = \left\{\begin{array}{ll}
        1, &  \text{if } \langle u, \xi^1_t(i)\rangle>0,\\
        0, &  \text{otherwise};
    \end{array}\right.~~~\,~~\\
    &\xi^2_t(i) = \left\{\begin{array}{ll}
    \bar{\xi}^2_t(i), &\text{ with probability 0.8,}\\
    1-\bar{\xi}^2_t(i), &\text{  with probability 0.2.}
    \end{array}\right.
\end{align*}
\end{subequations}
Clearly, for any $i\in\{1,2,\ldots,M\}$, $(\xi_t^1(i),\xi_t^2(i))_{t=1}^{\infty}$ forms a Markov chain. Let $\Pi_i$ denote the stationary distribution of the Markov chain on the $i$-th node. By defining the loss function as
$
\ell(x;\xi^1(i),\xi^2(i))=-\xi^2(i)\log(\sigma(\langle x, \xi^1(i) \rangle))-(1-\xi^2(i))\log(1-\sigma(\langle x,\xi^1(i)\rangle))
$
with $\sigma(t)=\frac{1}{1+\exp(-t)}$,
we  reconstruct $u$ as the solution to the following problem:
\begin{equation}
   \min_{x} ~~\sum\limits_{i=1}^m\EE_{(\xi^1(i),\xi^2(i))\sim\Pi_i}\ell(x;\xi^1(i),\xi^2(i)).
\end{equation}
 We choose $\gamma_k = \frac{1}{(k+1)^q}$ as our stepsize, where $q = 0.51$.  This choice is consistently with our theory below.
Specifically we compare:

{\bf DMGD}, where $(\xi^{1,k}(i),\xi^{2,k}(i))$ is from one trajectory of the Markov chain on the $i$-th node;

{\bf MCGD}, (i.e., the centralized Markov chain gradient descent), where   $(\xi^{1,k}(i),\xi^{2,k}(i))$ is from one trajectory of the Markov chain on the $i$-th node;

{\bf DSGD-$T$}, for $T=1,2,4,8,16$, where each $(\xi^{1,k}(i),\xi^{2,k}(i))$ is the $T$-th sample of an independent trajectory on the $i$-th node. All trajectories are generated by starting from the same initial state.

To compute $T$ gradients, DSGD-$T$ uses $T$ times as many samples as DMGD. We did not try to adapt $T$ as $k$ increases because there lacks a theoretical guidance. The numerical comparison results are reported in Figure 1, which show that DMGD outperforms the DSGD-$T$ with $T=1,2,4,8,16$.
The numerical results in Figure 1 are quite positive on DMGD. As  expected, DMGD used significantly fewer total samples than DSGD on each T. Surprisingly,  DMGD did not cost even more gradient computations.  It is important to find that DSGD1 and DSGD2, as well as DSGD4,  stagnate at noticeably lower accuracies due to that their T values are too small.  On the other hand, we observe that DMGD may not beat the centralized in the gradient computation and sampling costs, but it can significantly reduce the busiest node's communications.
\section{Proofs}
This section contains the proofs of the main results.

\subsection{Technical lemmas}
\begin{lemma}\label{lemc}
Let Assumption \ref{ass:mc} hold and let $\lambda_i(H)\in \mathbb{C}$ be  the $i$th largest eigenvalue of $H$, and
$$\lambda(H):=\frac{\max\{|\lambda_2(H)|,|\lambda_M(H)|\}+1}{2}\in [0,1).$$
Then, we can bound the largest entry-wise absolute value of the deviation matrix $\delta^k :=\Pi^*-H^{k}\in \mathbb{R}^{M\times M}$ as
\begin{equation}\label{core1}
    \|\delta^k\|_{\infty}\leq C_H\cdot\lambda^k(H)
\end{equation}
for $k\geq K$,  where $C_H$ is a constant that also depends on the Jordan canonical form of $H$ and  $K$ is a constant that depends on $\lambda(H)$ and $\lambda_2(H)$.
\end{lemma}
First, we define
${\bf 1}:=[1 1 \cdots 1 ]^{\top}$.
And the projection matrix is given as
\begin{equation}\label{prj}
    P:=\frac{\textbf{1}\textbf{1}^{\top}}{m}\in \mathbb{R}^{m\times m}.
\end{equation}
It is easy to see
\begin{equation}\label{pfacts}
    W P= P W=P.
\end{equation}
\begin{lemma}[Corollary 1.14., \cite{montenegro2006mathematical}]\label{estimate}
Let
$
    {\bf P}\in \mathbb{R}^{m\times m}
$ be the matrix whose elements are all $1/m$.
Given any $k\in \mathbb{Z}^+$, the mixing matrix ${\bf W}\in\RR^{m\times m}$ satisfies
$$\| W^k- P\|\leq \lambda_2(W)^k,$$
where $0\leq \lambda_2(W)<1$ is  the absolute value of  the second eigenvalue of $W$.
\end{lemma}

\begin{lemma}\label{bound1}
Let $({\bf x}^k)_{k\geq 0}$ be generated by   Algorithm 1 and Assumption \ref{ass1} hold,  then we have
\begin{equation}\label{bound1-r1}
 \frac{1}{m}\sum_{i=1}^m\|{\bf x}^{k}(i)- \overline{x^{k}}\|\leq \sqrt{m}B\sum_{j=0}^{k} \gamma_j\lambda_2(W)^{k-j}.
\end{equation}
for any $k\geq 0$. Further if $\gamma_j=\frac{1}{(j+1)^{\theta}}$ with $\frac{1}{2}\leq \theta<1$, it follows
\begin{equation}\label{bound1-r2}
 \frac{1}{m}\sum_{i=1}^m\|{\bf x}^{k}(i)- \overline{x^{k}}\|\leq \frac{\sqrt{m}BC_W}{(k+1)^{\theta}}
\end{equation}
and $C_W$ is a positive constant dependent on $W$ and $C_W=\mathcal{O}(\frac{1}{1-\lambda_2(W)})$.
\end{lemma}

\begin{lemma}[\cite{sun2018markov}]\label{lemcon}
Consider two nonnegative sequences $(\beta_k)_{k\geq 0}$ and $(h_k)_{k\geq 0}$ that satisfy
\begin{enumerate}
    \item $\lim_k h_k=0$ and $\sum_{k}h_k=+\infty$, and
    \item $\sum_{k}\beta_k h_k<+\infty$, and
    \item $|\beta_{k+1}-\beta_k|\leq c h_k$ for some $c>0$ and $k=0,1,\ldots$.
\end{enumerate}
Then, we have $\lim\beta_k=0$.
\end{lemma}

\begin{lemma}[\cite{sun2018markov}]\label{app}
Let $a>0$, and $x>0$ be a enough large real number. If
\begin{equation}\label{app-c}
    y-a\ln y+c=x.
\end{equation}
Then, it holds
\begin{equation}\label{app-r1}
    y-x\leq 2a\ln x.
\end{equation}
\end{lemma}

\subsection{Proof of  Lemma \ref{bound1}}
With direct calculating,
\begin{align*}
    &{\bf x}^k=W^k{\bf x}^0- \sum_{j=0}^{k-1} \gamma_j W^{k-j}{\bf u}^j.
\end{align*}
Recall \eqref{pfacts}, we have
\begin{align*}
    (I-W){\bf x}^k=(I-W)(I-P){\bf x}^k.
\end{align*}
Thus, we need to bound $\|(I-P){\bf x}^k\|$ and use the fact
$$(W- P)(I-P)= W-2 P+ W P= W-  P= W-  PW.$$
Using \eqref{pfacts}, we have
\begin{align*}
    &\|(I-P){\bf x}^k\|=\|(W^k-P){\bf x}^0- \sum_{j=0}^{k-1} \gamma_j (W^{k-j}-P){\bf u}^j\|\\
    &\leq \|W^k-P\|\cdot\|{\bf x}^0\|+\sum_{j=0}^{k-1} \gamma_j\| W^{k-j}-P\|\cdot\|{\bf u}^j\|\\
    &\leq B\sum_{j=0}^{k} \gamma_j\lambda_2(W)^{k-j}.
\end{align*}
Noticing that $\|(I-P){\bf x}^k\|^2=\frac{\sum_{i=1}^m \|{\bf x}^{k}(i)- \overline{x^{k}}\|^2}{m^2}$,
\begin{align*}
 &\frac{1}{m}\sum_{i=1}^m\|{\bf x}^{k}(i)- \overline{x^{k}}\|\leq\frac{\sqrt{m\sum_{i=1}^m \|{\bf x}^{k}(i)- \overline{x^{k}}\|^2}}{m}\\
 &\leq\sqrt{m}\|(I-P){\bf x}^k\|\leq \sqrt{m}B\sum_{j=0}^{k} \gamma_j\lambda_2(W)^{k-j}.
\end{align*}

\medskip

If $\gamma_k=\frac{1}{(k+1)^{\theta}}$, we have
\begin{align*}
&\sum_{j=0}^{k} \gamma_j\lambda_2(W)^{k-j}=\sum_{j=0}^{\ulcorner\frac{k}{2}\urcorner} \gamma_j\lambda_2(W)^{k-j}\\
&+\sum_{j=\ulcorner\frac{k}{2}\urcorner+1}^{k} \gamma_j\lambda_2(W)^{k-j}\leq k\lambda_2(W)^{k/2}\\
&+\frac{2^{\theta}}{(k+1)^{\theta}}\sum_{j=\ulcorner\frac{k}{2}\urcorner+1}^{k}\lambda_2(W)^{k-j}\leq k\lambda_2(W)^{k/2}\\
&+\frac{2^{\theta}}{(k+1)^{\theta}}\frac{1}{1-\lambda_2(W)}\leq \frac{C_W}{(k+1)^{\theta}},
\end{align*}
where $C_W:= \frac{2^{\theta}}{1-\lambda_2(W)}+\sup_{k}\{(k+1)^{1+\theta}\lambda_2(W)^{k/2}\}$ and $C_W=\mathcal{O}(\frac{1}{1-\lambda_2(W)})$.

\subsection{Proof of Theorem \ref{th2}}
Multiplying both sides of \eqref{dmgd2} with $P$,
 $$P{\bf x}^{k+1}= PW{\bf x}^k-\gamma_k P{\bf u}^k=P{\bf x}^k-\gamma_k P{\bf u}^k.$$
With direct computations, we get
\begin{align}\label{pronon}
    \| \overline{x^{k+1}}- \overline{x^{k}}\|=\frac{\|P{\bf x}^{k+1}-P{\bf x}^{k}\|}{\sqrt{m}}=\frac{\|\gamma_k P{\bf u}^k\|}{\sqrt{m}}\leq B\cdot \gamma_k .
\end{align}
 For integer $k\geq 1$, denote the  integer $\mathcal{T}_{k}$ as
\begin{align}\label{jianduanD}
 &\mathcal{T}_{k}:=\nonumber\\
 &\min\{\max\Big\{\Big\lceil\ln\Big(\frac{k}{ 2C_H B^2}\Big)/\ln(\frac{1}{\lambda(H)})\Big\rceil,K_H\Big\},k\}.
\end{align}
By using Lemma  \ref{lemc}, we then get
\begin{equation}\label{noncongd-p}
    \Big| [H^{\mathcal{T}_{k}}]_{i,j}-\frac{1}{M}\Big|\leq\frac{1/k}{2B^2},
\end{equation}
for any $i,j\in\{1,2,\ldots,M\}$.
The remaining of the proof consists of two major steps:
\begin{enumerate}
    \item in first step, we will prove $    \sum_{k}\gamma_k\EE\|\nabla f( \overline{x^{k-\mathcal{T}_{k}}})\|^2=O(\max\{1,\frac{1}{\ln(1/\lambda(H))}\})$, and
    \item in second step, we will show $\sum_{k}\big(\gamma_k\EE\|\nabla f( \overline{x^{k}})\|^2-\gamma_k\EE\|\nabla f(\overline{ x^{k-\mathcal{T}_{k}}})\|^2\big)=O(\max\{1,\frac{1}{\ln(1/\lambda(H))}\})$.
\end{enumerate}
Summing them together, we derive
\begin{equation}\label{th-sgd-t12}
    \sum_{k}\gamma_k\EE\|\nabla f( \overline{x^{k}})\|^2=\mathcal{O}(\frac{1}{\ln(1/\lambda(H))}\cdot \frac{1}{1-\lambda_2(W)}).
\end{equation}
Then, we   are led to
\begin{align}\label{th-sgd-t12-}
&(\sum_{i=1}^k\gamma_i)\cdot\EE(\min_{1\leq i\leq k}\{\|\nabla f( \overline{x^{i}})\|^2\})\leq \sum_{i=1}^k \gamma_{i} \EE\|\nabla f( \overline{x^i})\|^2\nonumber\\
&=\mathcal{O}(\frac{1}{\ln(1/\lambda(H))}\cdot \frac{1}{1-\lambda_2(W)}).
\end{align}
 Rearrangement of \eqref{th-sgd-t12-}  together with Schwarz inequality then gives us \eqref{th2-r2}.

\textbf{Step 1.}
Denote the shorthand notation
\begin{align*}
\small
\tilde{ u}^{k}:=\frac{\sum_{h=1}^m \nabla f_{j_{h,k}}({\bf x}^{k}(h)) }{m},~\tilde{ u}^{k-\mathcal{T}_{k}}:=\frac{\sum_{h=1}^m \nabla f_{j_{h,k}}({\bf x}^{k-\mathcal{T}_{k}}(h))}{m}
,
\end{align*}
we calculate the lower bound for following inner product:
\begin{align}\label{th-sgd-t3}
    &\EE_{\mathbf{j}_k}(\langle\nabla f( \overline{x^{k-\mathcal{T}_{k}}}),\tilde{ u}^{k-\mathcal{T}_{k}}\rangle\mid\chi^{k-\mathcal{T}_{k}})\nonumber\\
    &\quad\overset{a)}{=}\left\langle \nabla f( \overline{x^{k-\mathcal{T}_{k}}}), \right.\nonumber\\
    &\quad\quad\left.\frac{1}{m}\sum_{h=1}^m \sum_{i=1}^M \nabla f_{j_{h,k}}({\bf x}^{k-\mathcal{T}_{k}}(h))\cdot\mathbb{P}(j_{h,k}=i\mid\chi^{k-\mathcal{T}_{k}} )
    \right\rangle\nonumber\\
     &\quad\overset{b)}{=}\left\langle\nabla f( \overline{x^{k-\mathcal{T}_{k}}}), \right.\nonumber\\
&\quad\quad\left.\frac{1}{m}\sum_{h=1}^m \sum_{i=1}^M \nabla f_{j_{h,k}}({\bf x}^{k-\mathcal{T}_{k}}(h))\cdot\mathbb{P}(j_{h,k}=i\mid j_{h,k-\mathcal{T}_k} )
    \right\rangle\nonumber\\
     &\quad\overset{c)}{=}\left\langle\nabla f( \overline{x^{k-\mathcal{T}_{k}}}),\right.\nonumber\\
&\quad\quad\left.\frac{1}{m}\sum_{h=1}^m \sum_{i=1}^M \nabla f_{j_{h,k}}({\bf x}^{k-\mathcal{T}_{k}}(h))\cdot[H^{\mathcal{T}_{k}}]_{j_{m,k-\mathcal{T}_{k}},i}
    \right\rangle\nonumber\\
   &\quad\overset{d)}{\geq} \left\langle\nabla f( \overline{x^{k-\mathcal{T}_{k}}}), \frac{1}{m}\sum_{i=1}^m  \nabla f_{i}({\bf x}^{k-\mathcal{T}_{k}}(i))\right\rangle-\frac{1}{2k},
\end{align}
where $a)$ is due to the conditional expectation, and $b)$ uses the property of Markov chain, and $c)$ is  the matrix form of the probability, and $d)$ is from \eqref{noncongd-p}.
Direct calculations yield
\begin{align}\label{th-sgd-tt1'}
&\|  \nabla f( \overline{x^{k-\mathcal{T}_{k}}})\|^2=
   \left\langle\nabla f( \overline{x^{k-\mathcal{T}_{k}}}), \frac{1}{m}\sum_{i=1}^m  \nabla f_{i}({\bf x}^{k-\mathcal{T}_{k}}(i))\right\rangle\nonumber\\
   &=\|\nabla f( \overline{x^{k-\mathcal{T}_{k}}})\|^2+\langle\nabla f( \overline{x^{k-\mathcal{T}_{k}}}),\frac{1}{m}\sum_{i=1}^m  \nabla f_{i}({\bf x}^{k-\mathcal{T}_{k}}(i))\nonumber\\
    &-\nabla f( \overline{x^{k-\mathcal{T}_{k}}})\rangle\leq\|\nabla f( \overline{x^{k-\mathcal{T}_{k}}})\|^2\nonumber\\
    &+B\frac{L}{m}\sum_{i=1}^m \|{\bf x}^{k-\mathcal{T}_{k}}(i)- \overline{x^{k-\mathcal{T}_{k}}}\|.
\end{align}
Rearrangement of \eqref{th-sgd-t3} together with \eqref{th-sgd-tt1'} gives us
\begin{align}\label{th-sgd-tt1}
    &\gamma_k\EE\|  \nabla f( \overline{x^{k-\mathcal{T}_{k}}})\|^2\leq\gamma_k\EE(\langle  \nabla f( \overline{x^{k-\mathcal{T}_{k}}}),\tilde{ u}^{k-\mathcal{T}_{k}}\rangle)\nonumber\\
    &+B\frac{L}{m}\sum_{i=1}^m \gamma_k\EE\|{\bf x}^{k-\mathcal{T}_{k}}(i)- \overline{x^{k-\mathcal{T}_{k}}}\|+\frac{\gamma_k}{2k}.
\end{align}
We offer the bound for $ f( \overline{x^{k+1}})- f(\overline{x^{k}})$ as
\begin{align}\label{th-sgd-t4}
    &f( \overline{x^{k+1}})- f( x\overline{^{k}})\overset{a)}{\leq}\langle\nabla f( \overline{x^{k}}),  \overline{x^{k+1}}- \overline{x^{k}}\rangle\nonumber\\
    &+\frac{L}{2}\| \overline{x^{k+1}}-\overline{x^{k}}\|^2 \overset{b)}{=}\langle\nabla f( \overline{x^{k-\mathcal{T}_{k}}}), \overline{x^{k+1}}-\overline{x^{k}}\rangle\nonumber\\
    &+\langle\nabla f( \overline{x^{k}})-\nabla f( \overline{x^{k-\mathcal{T}_{k}}}), \overline{x^{k+1}}- \overline{x^{k}}\rangle+\frac{L}{2}\|\overline{x^{k+1}}-\overline{x^{k}}\|^2\nonumber\\
    &\overset{c)}{\leq}\langle\nabla f( \overline{x^{k-\mathcal{T}_{k}}}), \overline{x^{k+1}}- \overline{x^{k}}\rangle+\frac{(L+1)\| \overline{x^{k+1}}- \overline{x^{k}}\|^2}{2}\nonumber\\
    &+\frac{L^2\| \overline{x^k}- \overline{x^{k-\mathcal{T}_{k}}}\|^2}{2},
\end{align}
where $a)$ depends on the continuity of $\nabla f$, and $b)$ is the basic algebra computation, $c)$ uses the Schwarz inequality for $\langle\nabla f( \overline{x^{k}})-\nabla f( \overline{x^{k-\mathcal{T}_{k}}}), \overline{x^{k+1}}- \overline{x^{k}}\rangle$.
Moving $\langle\nabla f( \overline{x^{k-\mathcal{T}_{k}}}), \overline{x^{k+1}}- \overline{x^{k}}\rangle$ to left side,
\begin{align}\label{th-sgd-t6}
    &\langle\nabla f( \overline{x^{k-\mathcal{T}_{k}})}, \overline{x^k}- \overline{x^{k+1}}\rangle\leq  f( \overline{x^k})-f( \overline{x^{k+1}})\nonumber\\
    &\quad\quad+\frac{L^2\| \overline{x^{k+1}}- \overline{x^{k-\mathcal{T}_{k}}}\|^2}{2}+\frac{(L+1)\|\overline{x^{k+1}}- \overline{x^k}\|^2}{2}.
\end{align}
We then consider  the following bound:
\begin{align}\label{th-sgd-t7}
    &\EE(\langle\nabla f( \overline{x^{k-\mathcal{T}_{k}}}), \overline{x^k}- \overline{x^{k+1}}\rangle\mid\chi^{k-\mathcal{T}_{k}})\nonumber\\
    &\quad=\gamma_k\EE(\langle\nabla f( \overline{x^{k-\mathcal{T}_{k}}}), \tilde{v}^k\rangle\mid\chi^{k-\mathcal{T}_{k}})\nonumber\\
    &\quad=\gamma_k\EE(\langle \nabla f({\bf x}^{k-\mathcal{T}_{k}}), \tilde{ v}^{k-\mathcal{T}_{k}}\rangle\mid\chi^{k-\mathcal{T}_{k}})\nonumber\\
    &\quad+\gamma_k\EE(\langle\nabla f({\bf x}^{k-\mathcal{T}_{k}}),\tilde{v}^{k}-\tilde{ v}^{k-\mathcal{T}_{k}}\rangle\mid\chi^{k-\mathcal{T}_{k}})\nonumber\\
    &\quad\geq \gamma_k\EE(\langle \nabla f({\bf x}^{k-\mathcal{T}_{k}}),\tilde{ v}^{k-\mathcal{T}_{k}}\rangle\mid\chi^{k-\mathcal{T}_{k}})\nonumber\\
    &\quad-B\cdot L\cdot\EE(\gamma_k\| \overline{x^{k}}- \overline{x^{k-\mathcal{T}_{k}}}\|\mid\chi^{k-\mathcal{T}_{k}}),
\end{align}
where we used  the Lipschitz continuity and boundedness of $\nabla f$.
Taking conditional expectations on both sides of \eqref{th-sgd-t6} on $\chi^{k-\mathcal{T}_{k}}$ and rearrangement of \eqref{th-sgd-t7}, then we have
\begin{align}\label{th-sgd-t8}
   &\gamma_k\EE_{\mathbf{j}_k}(\langle\nabla f( \overline{x^{k-\mathcal{T}_{k}}}),\tilde{ v}^{k-\mathcal{T}_{k}}\rangle\mid\chi^{k-\mathcal{T}_{k}})\nonumber\\
   &\quad\leq\EE\big( f( \overline{x^k})-f( \overline{x^{k+1}})\mid\chi^{k-\mathcal{T}_{k}}\big)\nonumber\\
   &\quad+\frac{(L+1)\cdot\EE(\| \overline{x^{k+1}}-\overline{x^k}\|^2\mid\chi^{k-\mathcal{T}_{k}})}{2}\nonumber\\
   &\quad+B\cdot L\cdot\EE(\gamma_k\|\overline{x^{k}}-\overline{x^{k-\mathcal{T}_{k}}}\|\mid\chi^{k-\mathcal{T}_{k}})\nonumber\\
   &\quad+\frac{L^2\cdot\EE(\| \overline{x^{k+1}}- \overline{x^{k-\mathcal{T}_{k}}}\|^2\mid\chi^{k-\mathcal{T}_{k}})}{2}.
\end{align}
Taking expectations on both sides of  \eqref{th-sgd-t8} and with (\ref{th-sgd-tt1}), we are then led to
\begin{align}\label{th-sgd-t9}
   &\gamma_k\EE\|  \nabla f( \overline{x^{k-\mathcal{T}_{k}}})\|^2\leq\underbrace{\EE\big(f( \overline{x^k})-f(\overline{x^{k+1}})\big)}_{\textrm{(I)}}\nonumber\\
   &+\underbrace{\frac{(L+1)\cdot\EE\| \overline{x^{k+1}}-\overline{x^k}\|^2}{2}}_{\textrm{(II)}}+\underbrace{BL\cdot \gamma_k\cdot\EE\| \overline{x^{k}}-\overline{x^{k-\mathcal{T}_{k}}}\|}_{\textrm{(III)}}\nonumber\\
   &\quad+\underbrace{\frac{L^2\cdot\EE\|\overline{ x^{k+1}}-\overline{x^{k-\mathcal{T}_{k}}}\|^2}{2}}_{\textrm{(IV)}}\nonumber\\
   &\quad+\underbrace{B\frac{L}{m}\sum_{i=1}^m \gamma_k\EE\|{\bf x}^{k-\mathcal{T}_{k}}(i)- \overline{x^{k-\mathcal{T}_{k}}}\|}_{\textrm{(V)}}+\frac{\gamma_k}{2k}.
\end{align}
We now prove that $\textrm{(I)},\textrm{(II)},\textrm{(III)}$ and $\textrm{(IV)}$ are all summable.  The summability  $\textrm{(I)}$ is obvious. For $\textrm{(II)}$, $\textrm{(III)}$ and $\textrm{(IV)}$, with   \eqref{pronon}, we can derive (we omit the constant parameters in following)
\begin{align*}
    \textrm{(II)}: \EE\| \overline{x^{k+1}}- \overline{x^k}\|^2\leq \gamma_k B,
\end{align*}
and
\begin{align*}
\textrm{(III)}: &\EE(\gamma_k\| \overline{x^{k}}- \overline{x^{k-\mathcal{T}_{k}}}\|)\leq\gamma_k\sum_{d=k-\mathcal{T}_{k}}^{k-1}\EE\| \overline{x^{d+1}}- \overline{x^d}\|\nonumber\\
&\leq B\sum_{d=k-\mathcal{T}_{k}}^{k-1}\gamma_d\gamma_k\leq\frac{B}{2}\sum_{d=k-\mathcal{T}_{k}}^{k-1}(\gamma_d^2+\gamma_k^2)\nonumber\\
&=\frac{\mathcal{T}_{k} B}{2}\gamma_k^2+\frac{B}{2}\sum_{d=k-\mathcal{T}_{k}}^{k-1}\gamma_d^2,
\end{align*}
and
\begin{align*}
&\textrm{(IV)}: \EE(\| \overline{x^{k+1}}- \overline{x^{k-\mathcal{T}_{k}}}\|^2)\leq(\mathcal{T}_{k}+1)\sum_{d=k-\mathcal{T}_{k}}^{k}\EE\| \overline{x^{d+1}}- \overline{x^d}\|^2\nonumber\\
&\leq B^2(\mathcal{T}_{k}+1)\sum_{d=k-\mathcal{T}_{k}}^{k}\gamma_d^2
\end{align*}
and
\begin{align*}
&\textrm{(V)}: B\frac{L}{m}\sum_{i=1}^m \gamma_k\EE\|{\bf x}^{k-\mathcal{T}_{k}}(i)- \overline{x^{k-\mathcal{T}_{k}}}\|\\
&\leq B\frac{LC_W}{m}\gamma_{k-\mathcal{T}_{k}}\gamma_k\leq B\frac{LC_W}{2m}(\gamma_{k-\mathcal{T}_{k}}^2+\gamma_k^2)\\
&\leq B\frac{LC_W}{2m}\sum_{d=k-\mathcal{T}_{k}}^{k}\gamma_d^2.
\end{align*}
It is easy to see if $(\mathcal{T}_{k}\sum_{d=k-\mathcal{T}_{k}}^{k}\gamma_d^2)_{k\geq 0}$ is summable, $\textrm{(II)}$, $\textrm{(III)}$ and $\textrm{(IV)}$ are all summable. In fact, the proof of summability of $(\mathcal{T}_{k}\sum_{d=k-\mathcal{T}_{k}}^{k}\gamma_d^2)_{k\geq 0}$ can direct follow the proofs in \cite{sun2018markov}, which is presented here for completeness.

We consider a large enough integer $K$ which lets Lemma \ref{lemc} and Lemma \ref{app} be active, and $\mathcal{T}_{k}=\ln\Big(\frac{k}{ 2C_H B^2}\Big)/\ln(\frac{1}{\lambda(H)})$ when $k\geq K$.
Noting that  the summability of sequence is free of finite items, we then consider studying $(\mathcal{T}_{k}\sum_{d=k-\mathcal{J}_{k}}^{k-1}\gamma_d^2)_{k\geq K}$.
For any fixed integer $t\geq K$, $\gamma_t^2$ only  appears at index $k\geq K$ satisfying
$S_{t}:=\{k\in \mathbb{Z}^+\mid k-\mathcal{T}_{k}\leq t\leq k-1,~k\geq K\}$  in  the inner summation.
Let $k(t)$ be the solution of $k-\mathcal{T}_{k}=t$.
The direct computation gives us
$\sharp (S_t)\leq k(t)-t\leq 2\frac{\ln t}{\ln(1/\lambda(H))},$
where the last inequality comes from Lemma \ref{app}.
If $K$ is large enough, $\mathcal{T}_{k}\leq\frac{k}{2}$, and then
$k\leq 2t,~~\forall k\in S_t.$
Noting that $\mathcal{T}_{k}$ increases respect to $k$, we then get
$\mathcal{T}_{k}\leq \mathcal{T}_{2t},~~\forall k\in S_t.$
That means in $\sum_{k=K}^{+\infty}(\mathcal{T}_{k}\sum_{d=k-\mathcal{T}_{k}}^{k-1}\gamma_d^2)$,  $\gamma_t^2$ appears at most
$\mathcal{T}_{2t}\cdot \sharp(S_t)=O(\ln^2 t).$
The direct calculations then give us
\begin{align}\label{th-sgd-t-core}
    &\sum_{k=K}^{+\infty}\left(\mathcal{T}_{k}\sum_{d=k-\mathcal{T}_{k}}^{k-1}\gamma_d^2\right)\nonumber\\
    &=O( \sum_{t=K}\frac{\ln^2 t\cdot\gamma_t^2}{\ln(1/\lambda(H))})=\mathcal{O}(\ln(1/\lambda(H))).
\end{align}
Turning back to \eqref{th-sgd-t9} and Lemma \ref{bound1}, we are then led to
\begin{equation}\label{th-sgd-t10}
    \sum_{k}\gamma_k\EE\|  \nabla f(\overline{x^{k-\mathcal{T}_{k}}})\|^2=\mathcal{O}(\frac{1}{\ln(1/\lambda(H))}\cdot \frac{1}{1-\lambda_2(W)}).
\end{equation}

\textbf{Step 2:} With Lipschitz of $\nabla f$, it holds that
\begin{align}\label{th-sgd-t11}
    &\gamma_k\|\nabla f( \overline{x^{k}})\|^2-\gamma_k\|\nabla f( \overline{x^{k-\mathcal{T}_{k}}})\|^2\nonumber\\
    &\leq\gamma_k\langle\nabla f( \overline{x^{k}})-\nabla f( \overline{x^{k-\mathcal{T}_{k}}}),\nabla f( \overline{x^{k}})+\nabla f(\overline{x^{k-\mathcal{T}_{k}}})\rangle\nonumber\\
    &\leq\gamma_k\|\nabla f( \overline{x^{k}})-\nabla f(\overline{x^{k-\mathcal{T}_{k}}})\|\cdot\|\nabla f( \overline{x^{k}})+\nabla f( \overline{x^{k-\mathcal{T}_{k}}})\|\nonumber\\
    &\leq 2BL\gamma_k\| \overline{x^k}- \overline{x^{k-\mathcal{T}_{k}}}\|\nonumber\\
    &\leq BL\gamma_k^2+BL\| \overline{x^k}- \overline{x^{k-\mathcal{T}_{k}}}\|^2.
\end{align}
We have proved $\sum_{k\geq K}^{+\infty}\EE\| \overline{x^{k+1}}- \overline{x^{k-\mathcal{T}_{k}}}\|^2=O(\frac{1}{\ln(1/\lambda(H))})$; it is same way to prove  $\sum_{k\geq K}^{+\infty} \EE\| \overline{x^{k}}- \overline{x^{k-\mathcal{T}_{k}}}\|^2=O(\frac{1}{\ln(1/\lambda(H))})$. Thus, we can get
\begin{align*}
&\sum_{k}\big(\gamma_k\EE\|\nabla f( \overline{x^{k}})\|^2-\gamma_k\EE\|\nabla f(\overline{x^{k-\mathcal{T}_{k}}})\|^2\big)\nonumber\\
&=\mathcal{O}(\frac{1}{\ln(1/\lambda(H))}).
\end{align*}
On the other hand, with \eqref{pronon}, we can get
\begin{align}\label{nonergo-t1'}
    \mid\|\nabla f(\overline{x^{k+1}})\|^2-\|\nabla f(\overline{x^{k}})\|^2\mid=\mathcal{O}(\gamma_k).
\end{align}
Thus, we derive
\begin{align}\label{nonergo-t1}
    &\mid\EE\|\nabla f(\overline{x^{k+1}})\|^2-\EE\|\nabla f(\overline{ x^{k}})\|^2\mid\nonumber\\
    &\leq\EE\mid\|\nabla f( \overline{x^{k+1}})\|^2-\|\nabla f(\overline{x^{k}})\|^2\mid=\mathcal{O}(\gamma_k).
\end{align}
With \eqref{th-sgd-t12}  and \eqref{nonergo-t1}, Lemma \ref{lemcon} then gives us \eqref{th2-r1}.

\subsection{Proof of Theorem \ref{thzero}}
We   denote that
$${\bf h}^k:=\begin{bmatrix}
    ({\bf h}^{1,k})^{\top} \\
    ({\bf h}^{2,k})^{\top}\\
    \vdots \\
    ({\bf h}^{m,k})^{\top}\\
\end{bmatrix},~\,~\mathcal{G}^k:=\sigma({\bf x}^0,{\bf x}^1,\ldots,{\bf x}^k,{\bf u}^0,{\bf u}^1,\ldots,{\bf u}^k).$$
The Assumption 2 indicates $f$ is Lipschitz continuous, which together with \cite[Theorem 3.1]{ghadimi2013stochastic} directly gives
$$\left\|\EE\left(\frac{n(f({\bf x}+\delta {\bf h})-f({\bf x}))}{\delta} \mid {\bf h}\right)
-\nabla f({\bf x})\right\|=\mathcal{O}(\delta\cdot n^{\frac{3}{2}}).$$
With direct computations, we have
\begin{equation}\label{expesta}
    \EE( \overline{x^{k+1}}\mid {\bf h}^k)= \overline{x^k}-\gamma_k \frac{1}{m}\sum_{i=1}^m  \nabla f_{i}({\bf x}^{k}(i))+ e^k,
\end{equation}
where
\begin{align}\label{expesta2}
\| e^k\|=\mathcal{O}(\delta_k\cdot n^{\frac{3}{2}}).
\end{align}
We also denote that
$$\tilde{ u}^{k-\mathcal{T}_{k}}:=\frac{\sum_{h=1}^m \nabla f_{j_{h,k}}({\bf x}^{k-\mathcal{T}_{k}}(h))}{m}.$$
Now, we are prepared to prove the theorem. This proof is very similar to the proof of Theorem \ref{th2}. The main difference is to modify  \eqref{th-sgd-t7} as
\begin{align}\label{inexcore}
    &\EE(\langle\nabla f( \overline{x^{k-\mathcal{T}_{k}})},\overline{x^k}- \overline{x^{k+1}}\rangle\mid \mathcal{G}^k)\nonumber\\
    &\quad\overset{a)}{=}\gamma_k\EE(\langle\nabla f( \overline{x^{k-\mathcal{T}_{k}}}), \tilde{u}^k- e^k\rangle\mid\chi^{k-\mathcal{T}_{k}})\nonumber\\
    &\quad=\gamma_k\EE(\langle\nabla f( \overline{x^{k-\mathcal{T}_{k}}}), \tilde{u}^{k-\mathcal{T}_{k}}\rangle\mid\chi^{k-\mathcal{T}_{k}})\nonumber\\
    &\quad+\gamma_k\EE(\langle\nabla f(\overline{x^{k-\mathcal{T}_{k}}}), -e^k\rangle\mid\chi^{k-\mathcal{T}_{k}})\nonumber\\
    &\quad+\gamma_k\EE(\langle\nabla f( \overline{x^{k-\mathcal{T}_{k}}}),\tilde{u}^{k}-\tilde{u}^{k-\mathcal{T}_{k}}\rangle\mid\chi^{k-\mathcal{T}_{k}})\nonumber\\
    &\quad\overset{b)}{\geq} \gamma_k\EE(\langle \nabla f(\overline{x^{k-\mathcal{T}_{k}}}),\tilde{u}^{k-\mathcal{T}_{k}}\rangle\mid\chi^{k-\mathcal{T}_{k}})-B\gamma_k\cdot\|e^k\|\nonumber\\
    &\quad-B\cdot L\cdot\EE(\gamma_k\| \overline{x^{k}}- \overline{x^{k-\mathcal{T}_{k}}}\|\mid\chi^{k-\mathcal{T}_{k}}),
\end{align}
where $a)$ is due to \eqref{expesta}, and $b)$ depends on the Cauchy-Schwarz inequality. Taking total expectations on \eqref{inexcore}, we have
\begin{align}
   &\gamma_k\EE(\langle \nabla f(\overline{x^{k-\mathcal{T}_{k}}}),\tilde{u}^{k-\mathcal{T}_{k}}\rangle)\nonumber\\
   &\leq \textrm{(I)}+\textrm{(II)}+\textrm{(III)}+\textrm{(IV)}+B\gamma_k\cdot\| e^k\|,
\end{align}
where \textrm{(I)}, \textrm{(II)}, \textrm{(III)} and \textrm{(IV)} are given by \eqref{th-sgd-t9}. By using \eqref{expesta2}, we can prove
\begin{align}
   &\gamma_k\EE(\langle \nabla f( \overline{x^{k-\mathcal{T}_{k}}}),\tilde{u}^{k-\mathcal{T}_{k}}\rangle)\nonumber\\
   &=\mathcal{O}\left( \frac{1}{\ln(1/\lambda(H))}\cdot\frac{1}{1-\lambda_2(W)}+n^{\frac{3}{2}}\right).
\end{align}
The following is the same as  the proof for Theorem \ref{th2}.

\subsection{Proof of Proposition \ref{pro-online}}
The proofs of Proposition \ref{pro-online} is similar to the proof of Theorem \ref{th2}.
We first present a lemma as follows.
\begin{lemma}\label{bound2}
Let $({\bf x}^k)_{k\geq 0}$ be generated by   scheme \eqref{localscheme-online} and Assumption \ref{assf} hold, if   the stepsizes are selected as \eqref{stepsizes},  then we have
$$\frac{1}{m}\sum_{i=1}^m\|{\bf x}^{k}(i)- \overline{x^{k}}\|=\mathcal{O}(\frac{1}{1-\lambda_2(W)}\cdot{\frac{1}{(k+1)^{\theta}}}).$$
\end{lemma}
Like previous methods, the proof also consists of two parts:
\begin{enumerate}
    \item in the first one, we   prove $    \sum_{k}\gamma_k\EE\|\nabla \mathcal{F}(\overline{x^{k-\mathcal{T}_{k}}})\|^2=\mathcal{O}(\frac{1}{\ln(1/\lambda(H))}\cdot \frac{1}{1-\lambda_2(W)})$, and
    \item in second one, we focus on proving  $\sum_{k}\big(\gamma_k\EE\|\nabla \mathcal{F}( \overline{x^{k}})\|^2-\gamma_k\EE\|\nabla \mathcal{F}( \overline{x^{k-\mathcal{T}_{k}}})\|^2\big)=\mathcal{O}(\frac{1}{\ln(1/\lambda(H))}\cdot \frac{1}{1-\lambda_2(W)})$.
\end{enumerate}

\textbf{Step 1.}
Assume that $C^i$ and $\lambda_i$ are the factors in \eqref{core2} for Markov chain in the $i$th node. Let $C:=\max\{C^1,C^2,\ldots,C^m\}$ and  $\lambda:=\max\{\lambda_1,\lambda_2,\ldots,\lambda_m\}$.
For integer $k\geq 1$, we consider  the  integer $\mathcal{H}_{k}$ as
\begin{equation}\label{onlinet}
    \mathcal{H}_{k}:=\min\{\Big\lceil\ln\Big(\frac{k}{ 2C\cdot B^2}\Big)/\ln(1/\lambda)\Big\rceil,k\}.
\end{equation}
It is easy to see $\mathcal{H}_{k}\leq k$.
With [Theorem 4.9, \cite{montenegro2006mathematical}], we have the following relation
\begin{equation}\label{online1-t4}
\int_{\Xi}|p^{s+\mathcal{H}_{k}}_s(\xi)-\pi(\xi)| d\mu(\xi)\leq \frac{1}{2\cdot B^2\cdot k}, \forall s\in\mathbb{Z}^+
\end{equation}
where $p^{s+\mathcal{H}_{k}}_s(\xi)$ denotes the transition  p.d.f. from $s$ to $s+\mathcal{H}_{k}$ with respect to $\xi$. The property of time-homogeneous of the Markov chain directly gives that $p^{s+\mathcal{H}_{k}}_s(\xi)=p^{\mathcal{H}_{k}}_0 (\xi)$.
Denote the shorthand notation
$\tilde{ d}^{k-\mathcal{H}_{k}}:=\frac{1}{m}\sum_{i=1}^m \nabla f({\bf x}^{k-\mathcal{H}_{k}}(i);\xi^k(i)),$
we calculate the lower bound for following inner product:
\begin{align}\label{pro-sgd-t3}
    &\EE_{\mathbf{j}_k}(\langle\nabla \mathcal{F}( \overline{x^{k-\mathcal{H}_{k}}}),\tilde{d}^{k-\mathcal{H}_{k}}\rangle\mid\chi^{k-\mathcal{H}_{k}})\nonumber\\
    &\quad\overset{a)}{=}\left\langle\nabla \mathcal{F}( \overline{x^{k-\mathcal{H}_{k}}}),\sum_{i=1}^m \nabla F({\bf x}^{k-\mathcal{H}_{k}}(i);\xi) p_{k-\mathcal{H}_k}^{k}(\xi)d\mu(\xi)\right\rangle\nonumber\\
     &\quad\overset{b)}{=}\left\langle\nabla \mathcal{F}_{\alpha}({\bf x}^{k-\mathcal{H}_{k}}),\sum_{i=1}^m \nabla F({\bf x}^{k-\mathcal{H}_{k}}(i);\xi) p_{0}^{\mathcal{H}_k}(\xi)d\mu(\xi)\right\rangle\nonumber\\
   &\quad\overset{c)}{\geq}\|\nabla \mathcal{F}({\bf x}^{k-\mathcal{H}_{k}})\|^2-\frac{1}{2k}+\Big\langle\nabla \mathcal{F}_{\alpha}({\bf x}^{k-\mathcal{H}_{k}}),\nonumber\\
   &\quad\frac{1}{m}\sum_{i=1}^m [\nabla F({\bf x}^{k-\mathcal{H}_{k}}(i);\xi)-\nabla \mathcal{F}({\bf x}^{k-\mathcal{H}_{k}})] p_{0}^{\mathcal{H}_k}(\xi)d\mu(\xi)\Big\rangle\nonumber\\
   &\quad\overset{d)}{\geq}\|\nabla \mathcal{F}({\bf x}^{k-\mathcal{H}_{k}})\|^2-\frac{1}{2k}-\frac{B}{m}\sum_{i=1}^m\EE\|{\bf x}^{k}(i)- \overline{x^{k}}\|
\end{align}
where $a)$ uses the conditional expectation, and $b)$ comes from the property of Markov chain, and $c)$   depends on \eqref{noncongd-p}, and $d)$ is due to the Lipschitz property.
Rearrangement of \eqref{pro-sgd-t3} gives us
\begin{align}\label{pro-sgd-tt1}
    \gamma_k&\EE\|  \nabla \mathcal{F}( \overline{x^{k-\mathcal{H}_{k}}})\|^2\leq\gamma_k\EE(\langle  \nabla \mathcal{F}(\overline{x^{k-\mathcal{H}_{k}}}),\tilde{ d}^{k-\mathcal{H}_{k}}\rangle)\nonumber\\
    &+\frac{\gamma_k}{2k}+\gamma_k\frac{B}{m}\sum_{i=1}^m\EE\|{\bf x}^{k}(i)- \overline{x^{k}}\|.
\end{align}
Then we  need to  bound $ \mathcal{F}( \overline{x^{k+1}})- \mathcal{F}( \overline{x^{k}})$. With  Assumption \ref{assf},  $\nabla \mathcal{F}(\cdot)$ is Lipschitz continuous, the rest part is almost identical to the one of previous proof and will not be repeated.

\textbf{Step 2.} This step is very similar to  the second step in the proofs of  Theorem \ref{th2} and will not be repeated, either.
\section{Conclusions}
In this paper, we proposed the decentralized Markov chain gradient descent (DMGD) algorithm, where the samples are taken
along a trajectory of Markov chain over the network. Our algorithms can be used when it is impossible or very expensive to sample directly from a distribution, or the distribution is even unknown, but sampling via a Markov chain is possible and relatively cheap.
The convergence analysis is proved in possibly nonconvex cases.

Building upon the current work, several promising future directions can be pursued. The first one is to extend DMGD to the asynchronous setting, which can further reduce the synchronization overhead.
The second one is to reduce the communications cost in DMGD by using quantization or sparsification techniques.
Designing Markov chain primal-dual algorithms is also worth investigating.


\end{document}